\documentclass[12pt,reqno]{amsart}

\usepackage{amssymb}

\catcode`\@=11

\thinlines
\newskip\Einheit \Einheit=0.6cm
\newcount\xcoord \newcount\ycoord
\newdimen\xdim \newdimen\ydim \newdimen\PfadD@cke \newdimen\Pfadd@cke
\def\PfadDicke#1{\PfadD@cke#1 \divide\PfadD@cke by2 \Pfadd@cke\PfadD@cke \multiply\PfadD@cke by2}
\long\def\LOOP#1\REPEAT{\def\BODY{#1}\ITERATE}
\def\ITERATE{\BODY \let\next\ITERATE \else\let\next\relax\fi \next}
\let\REPEAT=\fi
\def\Punkt{\hbox{\raise-2pt\hbox to0pt{\hss\scriptsize$\bullet$\hss}}}
\def\DuennPunkt(#1,#2){\unskip
  \raise#2 \Einheit\hbox to0pt{\hskip#1 \Einheit
          \raise-2.5pt\hbox to0pt{\hss\normalsize$\bullet$\hss}\hss}}
\def\NormalPunkt(#1,#2){\unskip
  \raise#2 \Einheit\hbox to0pt{\hskip#1 \Einheit
          \raise-3pt\hbox to0pt{\hss\large$\bullet$\hss}\hss}}
\def\DickPunkt(#1,#2){\unskip
  \raise#2 \Einheit\hbox to0pt{\hskip#1 \Einheit
          \raise-4pt\hbox to0pt{\hss\Large$\bullet$\hss}\hss}}
\def\Kreis(#1,#2){\unskip
  \raise#2 \Einheit\hbox to0pt{\hskip#1 \Einheit
          \raise-4pt\hbox to0pt{\hss\Large$\circ$\hss}\hss}}
\def\Diagonale(#1,#2)#3{\unskip\leavevmode
  \xcoord#1\relax \ycoord#2\relax
      \raise\ycoord \Einheit\hbox to0pt{\hskip\xcoord \Einheit
         \unitlength\Einheit
         \line(1,1){#3}\hss}}
\def\AntiDiagonale(#1,#2)#3{\unskip\leavevmode
  \xcoord#1\relax \ycoord#2\relax 
      \raise\ycoord \Einheit\hbox to0pt{\hskip\xcoord \Einheit
         \unitlength\Einheit
         \line(1,-1){#3}\hss}}
\def\Pfad(#1,#2),#3\endPfad{\unskip\leavevmode
  \xcoord#1 \ycoord#2 \thicklines\ZeichnePfad#3\endPfad\thinlines}
\def\ZeichnePfad#1{\ifx#1\endPfad\let\next\relax
  \else\let\next\ZeichnePfad
    \ifnum#1=1
      \raise\ycoord \Einheit\hbox to0pt{\hskip\xcoord \Einheit
         \vrule height\Pfadd@cke width1 \Einheit depth\Pfadd@cke\hss}%
      \advance\xcoord by 1
    \else\ifnum#1=2
      \raise\ycoord \Einheit\hbox to0pt{\hskip\xcoord \Einheit
        \hbox{\hskip-\PfadD@cke\vrule height1 \Einheit width\PfadD@cke depth0pt}\hss}%
      \advance\ycoord by 1
    \else\ifnum#1=3
      \raise\ycoord \Einheit\hbox to0pt{\hskip\xcoord \Einheit
         \unitlength\Einheit
         \line(1,1){1}\hss}
      \advance\xcoord by 1
      \advance\ycoord by 1
    \else\ifnum#1=4
      \raise\ycoord \Einheit\hbox to0pt{\hskip\xcoord \Einheit
         \unitlength\Einheit
         \line(1,-1){1}\hss}
      \advance\xcoord by 1
      \advance\ycoord by -1
    \else\ifnum#1=5
      \advance\xcoord by -1
      \raise\ycoord \Einheit\hbox to0pt{\hskip\xcoord \Einheit
         \vrule height\Pfadd@cke width1 \Einheit depth\Pfadd@cke\hss}%
    \else\ifnum#1=6
      \advance\ycoord by -1
      \raise\ycoord \Einheit\hbox to0pt{\hskip\xcoord \Einheit
        \hbox{\hskip-\PfadD@cke\vrule height1 \Einheit width\PfadD@cke depth0pt}\hss}%
    \else\ifnum#1=7
      \advance\xcoord by -1
      \advance\ycoord by -1
      \raise\ycoord \Einheit\hbox to0pt{\hskip\xcoord \Einheit
         \unitlength\Einheit
         \line(1,1){1}\hss}
    \else\ifnum#1=8
      \advance\xcoord by -1
      \advance\ycoord by +1
      \raise\ycoord \Einheit\hbox to0pt{\hskip\xcoord \Einheit
         \unitlength\Einheit
         \line(1,-1){1}\hss}
    \fi\fi\fi\fi
    \fi\fi\fi\fi
  \fi\next}
\def\hSSchritt{\leavevmode\raise-.4pt\hbox to0pt{\hss.\hss}\hskip.2\Einheit
  \raise-.4pt\hbox to0pt{\hss.\hss}\hskip.2\Einheit
  \raise-.4pt\hbox to0pt{\hss.\hss}\hskip.2\Einheit
  \raise-.4pt\hbox to0pt{\hss.\hss}\hskip.2\Einheit
  \raise-.4pt\hbox to0pt{\hss.\hss}\hskip.2\Einheit}
\def\vSSchritt{\vbox{\baselineskip.2\Einheit\lineskiplimit0pt
\hbox{.}\hbox{.}\hbox{.}\hbox{.}\hbox{.}}}
\def\DSSchritt{\leavevmode\raise-.4pt\hbox to0pt{%
  \hbox to0pt{\hss.\hss}\hskip.2\Einheit
  \raise.2\Einheit\hbox to0pt{\hss.\hss}\hskip.2\Einheit
  \raise.4\Einheit\hbox to0pt{\hss.\hss}\hskip.2\Einheit
  \raise.6\Einheit\hbox to0pt{\hss.\hss}\hskip.2\Einheit
  \raise.8\Einheit\hbox to0pt{\hss.\hss}\hss}}
\def\dSSchritt{\leavevmode\raise-.4pt\hbox to0pt{%
  \hbox to0pt{\hss.\hss}\hskip.2\Einheit
  \raise-.2\Einheit\hbox to0pt{\hss.\hss}\hskip.2\Einheit
  \raise-.4\Einheit\hbox to0pt{\hss.\hss}\hskip.2\Einheit
  \raise-.6\Einheit\hbox to0pt{\hss.\hss}\hskip.2\Einheit
  \raise-.8\Einheit\hbox to0pt{\hss.\hss}\hss}}
\def\SPfad(#1,#2),#3\endSPfad{\unskip\leavevmode
  \xcoord#1 \ycoord#2 \ZeichneSPfad#3\endSPfad}
\def\ZeichneSPfad#1{\ifx#1\endSPfad\let\next\relax
  \else\let\next\ZeichneSPfad
    \ifnum#1=1
      \raise\ycoord \Einheit\hbox to0pt{\hskip\xcoord \Einheit
         \hSSchritt\hss}%
      \advance\xcoord by 1
    \else\ifnum#1=2
      \raise\ycoord \Einheit\hbox to0pt{\hskip\xcoord \Einheit
        \hbox{\hskip-2pt \vSSchritt}\hss}%
      \advance\ycoord by 1
    \else\ifnum#1=3
      \raise\ycoord \Einheit\hbox to0pt{\hskip\xcoord \Einheit
         \DSSchritt\hss}
      \advance\xcoord by 1
      \advance\ycoord by 1
    \else\ifnum#1=4
      \raise\ycoord \Einheit\hbox to0pt{\hskip\xcoord \Einheit
         \dSSchritt\hss}
      \advance\xcoord by 1
      \advance\ycoord by -1
    \else\ifnum#1=5
      \advance\xcoord by -1
      \raise\ycoord \Einheit\hbox to0pt{\hskip\xcoord \Einheit
         \hSSchritt\hss}%
    \else\ifnum#1=6
      \advance\ycoord by -1
      \raise\ycoord \Einheit\hbox to0pt{\hskip\xcoord \Einheit
        \hbox{\hskip-2pt \vSSchritt}\hss}%
    \else\ifnum#1=7
      \advance\xcoord by -1
      \advance\ycoord by -1
      \raise\ycoord \Einheit\hbox to0pt{\hskip\xcoord \Einheit
         \DSSchritt\hss}
    \else\ifnum#1=8
      \advance\xcoord by -1
      \advance\ycoord by 1
      \raise\ycoord \Einheit\hbox to0pt{\hskip\xcoord \Einheit
         \dSSchritt\hss}
    \fi\fi\fi\fi
    \fi\fi\fi\fi
  \fi\next}
\def\Koordinatenachsen(#1,#2){\unskip
 \hbox to0pt{\hskip-.5pt\vrule height#2 \Einheit width.5pt depth1 \Einheit}%
 \hbox to0pt{\hskip-1 \Einheit \xcoord#1 \advance\xcoord by1
    \vrule height0.25pt width\xcoord \Einheit depth0.25pt\hss}}
\def\Koordinatenachsen(#1,#2)(#3,#4){\unskip
 \hbox to0pt{\hskip-.5pt \ycoord-#4 \advance\ycoord by1
    \vrule height#2 \Einheit width.5pt depth\ycoord \Einheit}%
 \hbox to0pt{\hskip-1 \Einheit \hskip#3\Einheit 
    \xcoord#1 \advance\xcoord by1 \advance\xcoord by-#3 
    \vrule height0.25pt width\xcoord \Einheit depth0.25pt\hss}}
\def\Gitter(#1,#2){\unskip \xcoord0 \ycoord0 \leavevmode
  \LOOP\ifnum\ycoord<#2
    \loop\ifnum\xcoord<#1
      \raise\ycoord \Einheit\hbox to0pt{\hskip\xcoord \Einheit\Punkt\hss}%
      \advance\xcoord by1
    \repeat
    \xcoord0
    \advance\ycoord by1
  \REPEAT}
\def\Gitter(#1,#2)(#3,#4){\unskip \xcoord#3 \ycoord#4 \leavevmode
  \LOOP\ifnum\ycoord<#2
    \loop\ifnum\xcoord<#1
      \raise\ycoord \Einheit\hbox to0pt{\hskip\xcoord \Einheit\Punkt\hss}%
      \advance\xcoord by1
    \repeat
    \xcoord#3
    \advance\ycoord by1
  \REPEAT}
\def\Label#1#2(#3,#4){\unskip \xdim#3 \Einheit \ydim#4 \Einheit
  \def\lo{\advance\xdim by-.5 \Einheit \advance\ydim by.5 \Einheit}%
  \def\llo{\advance\xdim by-.25cm \advance\ydim by.5 \Einheit}%
  \def\loo{\advance\xdim by-.5 \Einheit \advance\ydim by.25cm}%
  \def\o{\advance\ydim by.25cm}%
  \def\ro{\advance\xdim by.5 \Einheit \advance\ydim by.5 \Einheit}%
  \def\rro{\advance\xdim by.25cm \advance\ydim by.5 \Einheit}%
  \def\roo{\advance\xdim by.5 \Einheit \advance\ydim by.25cm}%
  \def\l{\advance\xdim by-.30cm}%
  \def\r{\advance\xdim by.30cm}%
  \def\lu{\advance\xdim by-.5 \Einheit \advance\ydim by-.6 \Einheit}%
  \def\llu{\advance\xdim by-.25cm \advance\ydim by-.6 \Einheit}%
  \def\luu{\advance\xdim by-.5 \Einheit \advance\ydim by-.30cm}%
  \def\u{\advance\ydim by-.30cm}%
  \def\ru{\advance\xdim by.5 \Einheit \advance\ydim by-.6 \Einheit}%
  \def\rru{\advance\xdim by.25cm \advance\ydim by-.6 \Einheit}%
  \def\ruu{\advance\xdim by.5 \Einheit \advance\ydim by-.30cm}%
  #1\raise\ydim\hbox to0pt{\hskip\xdim
     \vbox to0pt{\vss\hbox to0pt{\hss$#2$\hss}\vss}\hss}%
}
\catcode`\@=12

\font\HUGE=cmmi12 scaled 4000

\numberwithin{equation}{section}

\newtheorem{Theorem}{Theorem}

\newtheorem{Problem}{Problem}

\theoremstyle{definition}

\theoremstyle{remark}

\newtheorem*{Remarks}{Remarks}

\def\la{\lambda}
\def\rh{\rho}

\def\De{\Delta}

\def\today{\ifcase\month\or
 January\or February\or March\or April\or May\or June\or
 July\or August\or September\or October\or November\or December\fi
 \space\number\day, \number\year}

\def\({\left(}
\def\){\right)}
\def\[{\left[}
\def\]{\right]}

\def\3{\ss}

\def\cross{\operatorname{cross}}
\def\nest{\operatorname{nest}}
\def\ccross{\operatorname{\overline{\vphantom{t}cross}}}
\def\nnest{\operatorname{\overline{nest}}}
\def\NE{\text{\it NE}}
\def\SE{\text{\it SE}}
\def\Ne{\text{\it Ne}}
\def\Se{\text{\it Se}}

\begin{document}

\newbox\Adr
\setbox\Adr\vbox{
\centerline{\sc C.~Krattenthaler$^\dagger$}
\vskip18pt
\centerline{Fakult\"at f\"ur Mathematik, Universit\"at Wien,}
\centerline{Nordbergstra\3e 15, A-1090 Vienna, Austria.}
\centerline{\footnotesize WWW: \footnotesize\tt http://www.mat.univie.ac.at/\~{}kratt}
}

\title[Growth diagrams]{Growth diagrams, and 
increasing and decreasing
chains in fillings of Ferrers shapes}

\author[C.~Krattenthaler]{\box\Adr}

\address{Fakult\"at f\"ur Mathematik, Universit\"at Wien,
Nordbergstrasze 15, A-1090 Vienna, Austria.
WWW: \tt http://www.mat.univie.ac.at/\~{}kratt}

\dedicatory{Dedicated to Amitai Regev on the occasion of his 65th
birthday}

\thanks{$^\dagger$ Research partially supported by EC's IHRP Programme,
grant HPRN-CT-2001-00272, ``Algebraic Combinatorics in Europe", 
by the Austrian Science Foundation FWF, grant S9607-N13,
in the framework of the National Research Network
``Analytic Combinatorics and Probabilistic Number Theory,"
and by
the ``Algebraic Combinatorics" Programme during Spring 2005 
of the Institut Mittag--Leffler of the Royal Swedish Academy of Sciences}


\begin{abstract}
We put recent results by Chen, Deng, Du, Stanley and Yan on
crossings and nestings of matchings and set partitions in the
larger context of the enumeration of fillings of Ferrers shape
on which one imposes restrictions on their increasing and decreasing
chains. While Chen et al.\ work with Robinson--Schensted-like
insertion/deletion algorithms, we use the growth diagram construction
of Fomin to obtain our results. We extend the results by Chen et al.,
which, in the language of fillings, are results about
$0$-$1$-fillings, to arbitrary fillings. Finally, we point out that,
very likely, these results are part of a bigger picture which also
includes recent results of Jonsson on $0$-$1$-fillings of stack
polyominoes, and of results of Backelin, West and Xin and of
Bousquet--M\'elou and Steingr\'\i msson on the enumeration of
permutations and involutions with restricted patterns.
In particular, we show that our growth diagram bijections do
in fact provide alternative proofs of the results by
Backelin, West and Xin and by Bousquet--M\'elou and Steingr\'\i msson.
\end{abstract}

\maketitle

\section{Introduction} 
In the recent paper \cite{ChDDAB}, Chen et al.\ used
Robinson--Schensted like insertion/deletion processes to prove
enumeration results about matchings and, more generally, set
partitions with certain restrictions on their crossings and nestings.
At the heart of their results, there is Greene's
theorem \cite{GreCAA} on the relation between increasing and
decreasing subsequences in permutations and the shape of the
tableaux which are obtained by the Robinson--Schensted
correspondence.
The purpose of this paper is to put these results in a larger
context, namely the context of enumeration of fillings of shapes
where one imposes restrictions on the increasing and decreasing
chains of the fillings. 

As we explain in Section~\ref{sec:part}, 
the results in \cite{ChDDAB} are equivalent to results about
$0$-$1$-fillings of triangular arrangements of cells.
We show in Section~\ref{sec:growth}
that (almost) all the theorems of \cite{ChDDAB} generalize
to the setting of $0$-$1$-fillings of arrangements of cells
which have the form of a Ferrers diagram. 
In contrast to Chen et al., we work with the growth diagram construction of 
Fomin \cite{FomiAZ} which, in my opinion, leads to a more transparent 
presentation of the bijections between fillings and oscillating
sequences of (integer) partitions which underlie these results.

In Section~\ref{sec:RSK}, we extend these constructions to the
``Knuth setting," that is, to fillings where the constraint that in 
each row and in each column there is at most one $1$ is dropped.
The results that can be obtained in this way are presented in
Theorems~\ref{thm:1a}--\ref{thm:2a}. They are based on four
variations of the Robinson--Schensted--Knuth correspondence
(presented here in the language of growth diagrams). These
four variations are, in one form or another, well-known (see
\cite[App.~A.4]{FultAC}), and their growth diagram versions have been
worked out by Roby \cite{RobyAA} (for the original
Robinson--Schensted--Knuth correspondence) and van Leeuwen
\cite{LeeuAH}.
For the convenience of the reader, we recall these correspondences
(in growth diagram version and in insertion version) in
Section~\ref{sec:RSK}, and we use the opportunity to also
work out the corresponding variations of Greene's theorem on 
increasing and decreasing chains, see
Theorems~\ref{thm:3a}, \ref{thm:3b} and \ref{thm:3d}, which so far
has only been done for the original Robinson--Schensted--Knuth
correspondence.

The motivation to go beyond the setting of Chen et al.\ comes from
recent results by Jonsson \cite{JonsAA} on the enumeration
of $0$-$1$-fillings of stack polyominoes (shapes, which are more
general than Ferrers diagrams)
with restricted length of increasing and
decreasing chains of $1$'s. The restriction that in each row and in
each column of the shape there is at most one $1$ is not
present in the results by Jonsson. Unfortunately, the result on
$0$-$1$-fillings (see Theorem~\ref{thm:2a}, Eq.~\eqref{eq:NES2})
which we are able to obtain from the growth diagram construction
does not imply Jonsson's result, although it comes very close.
Nevertheless, we believe that all these phenomena are part of a
bigger picture that needs to be uncovered.
To this bigger picture there belong certainly results by
Backelin, West and Xin \cite{BaWXAA} and by
Bousquet--M\'elou and Steingr\'\i msson \cite{BoStAA} on the
enumeration of permutations and involutions with restricted
patterns. Indeed, we show that our growth diagram bijections provide
alternative proofs for their results. (I owe this observation to
Mireille Bousquet--M\'elou.) In spite of this,
apparently we do not yet have the right understanding for 
all these phenomena (not even regarding \cite{BaWXAA} and
\cite{BoStAA}). 
In particular, this paper seems to indicate that
Robinson--Schensted--Knuth-like insertion/deletion,
respectively growth diagram processes,
do not seem to be right tools for proving the results from
\cite{JonsAA}, say. 
All this is made more precise in Section~\ref{sec:big},
where we also formulate several open problems that should lead to
an understanding of the ``big picture."

To conclude the introduction, we point out that, independently from
\cite{ChDDAB}, Robinson--Schensted
like algorithms between set partitions and oscillating sequences of
(integer) partitions have also been constructed by
Halverson and Lewandowski \cite{HaLeAA}, with the completely
different motivation of explaining combinatorial identities 
arising from the representation theory of the partition algebra.
Halverson and Lewandowski provide both the insertion/deletion and the
growth diagram presentation of the algorithms. However, in their
considerations, Greene's theorem does not play any role.

\section{Growth diagrams} \label{sec:growth}

In this section we review the basic facts on growth diagrams
(cf.\ \cite[Sec.~3]{BrFoAA},
\cite[Ch.~2]{RobyAA}, \cite{RobyAD}, \cite[Sec.~5.2]{SagaAQ}, 
\cite[Sec.~7.13]{StanBI}). 

We start by fixing the standard partition notation (cf.\
e.g.\ \cite[Sec.~7.2]{StanBI}).
A {\it partition} is a weakly decreasing sequence
$\la=(\la_1,\la_2,\dots,\la_\ell)$ of positive integers.
This also includes the {\it empty partition} $()$, denoted by
$\emptyset$. For the sake of convenience, we shall often tacitly
identify a partition $\la=(\la_1,\la_2,\dots,\la_\ell)$ with the
infinite sequence $(\la_1,\la_2,\dots,\la_\ell,0,0,\dots)$, that is,
the sequence which arises from $\la$ by appending infinitely
many $0$'s. To each partition $\la$, one
associates its {\it Ferrers diagram} (also called {\it Ferrers shape}), 
which is the left-justified
arrangement of squares with $\la_i$ squares in the $i$-th row,
$i=1,2,\dots$. 
We define a {\it partial order} $\subseteq$ on partitions 
by containment of their Ferrers diagrams.
The {\it union} $\mu\cup\nu$ of two partitions $\mu$ and $\nu$
is the partition which arises by forming the union of the Ferrers
diagrams of $\mu$ and $\nu$. Thus, if
$\mu=(\mu_1,\mu_2,\dots)$ and
$\nu=(\nu_1,\nu_2,\dots)$, then $\mu\cup\nu$ is the partition 
$\la=(\la_1,\la_2,\dots)$, where $\la_i=\max\{\mu_i,\nu_i\}$ for 
$i=1,2,\dots$. The {\it intersection} 
$\mu\cap\nu$ of two partitions $\mu$ and $\nu$
is the partition which arises by forming the intersection of the Ferrers
diagrams of $\mu$ and $\nu$. Thus, if
$\mu=(\mu_1,\mu_2,\dots)$ and
$\nu=(\nu_1,\nu_2,\dots)$, then $\mu\cap\nu$ is the partition 
$\rh=(\rh_1,\rh_2,\dots)$, where $\rh_i=\min\{\mu_i,\nu_i\}$ for 
$i=1,2,\dots$.
The {\em conjugate of\/} a partition $\lambda$ is the partition 
$(\lambda^\prime_1, \dots, \lambda^\prime_{\lambda_1})$ where
$\lambda_j'$ is the length of the $j$-th column in the Ferrers
diagram of $\lambda$.

The objects of consideration in the present paper
are fillings of arrangements of cells which look like Ferrers
shapes in French notation, that is, which have straight left side,
straight bottom side, and which support a descending staircase, see
Figure~\ref{fig:00}.a for an example. We shall encode such Ferrers
shape by sequences of $D$'s and $R$'s, where $D$ stands for ``down"
and $R$ stands for ``right." To construct the sequence of $D$'s and
$R$'s of a Ferrers shape we trace the right/up boundary of the
Ferrers shape from top-left to bottom-right and write $D$ whenever we
encounter a down-step, respectively $R$ whenever we encounter a
right-step. For example, the Ferrers shape in Figure~\ref{fig:00}.a
would be encoded by $RDRDDRDDRRD$.

\begin{figure}[h]
$$
\Einheit.4cm
\Pfad(0,18),666666666666666666111111111111111\endPfad
\Pfad(0,18),111666666666666666666\endPfad
\Pfad(0,15),111111666666666666666\endPfad
\Pfad(0,12),111111666111666666666\endPfad
\Pfad(0,9),111111111666666111666\endPfad
\Pfad(0,6),111111111666111111666\endPfad
\Pfad(0,3),111111111111111666\endPfad
\hbox{\hskip7.5cm}
\Pfad(0,18),666666666666666666111111111111111\endPfad
\Pfad(0,18),111666666666666666666\endPfad
\Pfad(0,15),111111666666666666666\endPfad
\Pfad(0,12),111111666111666666666\endPfad
\Pfad(0,9),111111111666666111666\endPfad
\Pfad(0,6),111111111666111111666\endPfad
\Pfad(0,3),111111111111111666\endPfad
\Label\ro{\text {\Huge $0$}}(1,1)
\Label\ro{\text {\Huge $0$}}(4,1)
\Label\ro{\text {\Huge $0$}}(7,1)
\Label\ro{\text {\Huge $0$}}(10,1)
\Label\ro{\text {\Huge $0$}}(1,4)
\Label\ro{\text {\Huge $0$}}(7,4)
\Label\ro{\text {\Huge $0$}}(1,7)
\Label\ro{\text {\Huge $0$}}(4,7)
\Label\ro{\text {\Huge $0$}}(7,7)
\Label\ro{\text {\Huge $0$}}(4,10)
\Label\ro{\text {\Huge $0$}}(1,13)
\Label\ro{\text {\Huge $0$}}(4,13)
\Label\ro{\text {\Huge $0$}}(1,16)
\Label\ro{\text {\Huge $1$}}(4,4)
\Label\ro{\text {\Huge $1$}}(1,10)
\Label\ro{\text {\Huge $1$}}(13,1)
\hskip6cm
$$
\vskip10pt
\centerline{\small a. A Ferrers shape in French notation
\hskip1cm
b. A filling of the Ferrers shape}
\caption{}
\label{fig:00}
\end{figure}

We fill the cells of such a Ferrers shape $F$ with non-negative integers.
In this section (and in the following section) 
the fillings will be restricted to $0$-$1$-fillings
such that every row and every column contains at most one $1$.
See Figure~\ref{fig:00}.b for an example.

Next, the corners of the cells are labelled by partitions
such that the following two conditions are satisfied:

\begin{enumerate}
\item[(C1)] A partition is either equal to its right neighbour
or smaller by exactly one square, the same being true for a
partition and its top neighbour.
\item[(C2)] A partition and its right neighbour are equal if and only
if in the column of cells of $F$ below them there appears no $1$ and if
their bottom neighbours are also equal to each other.
Similarly, a partition and its top neighbour are equal if and only if
in the row of cells of $F$ to the left of them there appears no $1$ and if
their left neighbours are also equal to each other.
\end{enumerate}

\begin{figure}[h]
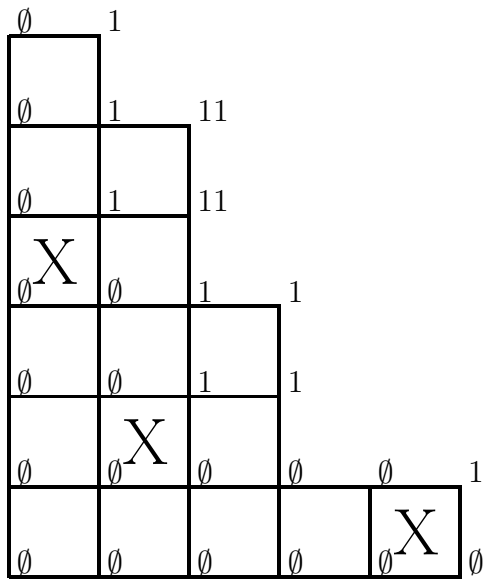

$$
\Einheit.4cm
\Pfad(0,18),666666666666666666111111111111111\endPfad
\Pfad(0,18),111666666666666666666\endPfad
\Pfad(0,15),111111666666666666666\endPfad
\Pfad(0,12),111111666111666666666\endPfad
\Pfad(0,9),111111111666666111666\endPfad
\Pfad(0,6),111111111666111111666\endPfad
\Pfad(0,3),111111111111111666\endPfad
\Label\ro{\emptyset}(0,0)
\Label\ro{\emptyset}(0,3)
\Label\ro{\emptyset}(0,6)
\Label\ro{\emptyset}(0,9)
\Label\ro{\emptyset}(0,12)
\Label\ro{\emptyset}(0,15)
\Label\ro{\emptyset}(0,18)
\Label\ro{1}(3,18)
\Label\ro{1}(3,15)
\Label\ro{1}(3,12)
\Label\ro{\emptyset}(3,9)
\Label\ro{\emptyset}(3,6)
\Label\ro{\emptyset}(3,3)
\Label\ro{\emptyset}(3,0)
\Label\ro{\hphantom{1}11}(6,15)
\Label\ro{\hphantom{1}11}(6,12)
\Label\ro{1}(6,9)
\Label\ro{1}(6,6)
\Label\ro{\emptyset}(6,3)
\Label\ro{\emptyset}(6,0)
\Label\ro{1}(9,9)
\Label\ro{1}(9,6)
\Label\ro{\emptyset}(9,3)
\Label\ro{\emptyset}(9,0)
\Label\ro{\emptyset}(12,3)
\Label\ro{\emptyset}(12,0)
\Label\ro{1}(15,3)
\Label\ro{\emptyset}(15,0)
\Label\ro{\text {\Huge X}}(4,4)
\Label\ro{\text {\Huge X}}(1,10)
\Label\ro{\text {\Huge X}}(13,1)
\hskip6cm
$$
\caption{A growth diagram}
\label{fig:0}
\end{figure}

See Figure~\ref{fig:0} for an example. (More examples can be found in
Figures~\ref{fig:3}--\ref{fig:5}.)
There, we use a short notation for partitions. For example,
$11$ is short for $(1,1)$.
Moreover, we changed the convention of representing the filling slightly
for better visibility,
by suppressing $0$'s and by replacing $1$'s by X's. Indeed, the
filling represented in Figure~\ref{fig:0} is the same as the one in
Figure~\ref{fig:00}.b.

Diagrams which obey the conditions (C1) and (C2) are called {\it
growth diagrams}.

We are interested in growth diagrams which
obey the following ({\it forward}) {\it local rules} 
(see Figure~\ref{fig:1}).

\begin{figure}[h]
$$
\Pfad(0,0),1111222255556666\endPfad
\Label\lu{\rh}(0,0)
\Label\ru{\mu}(4,0)
\Label\lo{\nu}(0,4)
\Label\ro{\la}(4,4)
\hbox{\hskip6cm}
\Pfad(0,0),1111222255556666\endPfad
\Label\lu{\rh}(0,0)
\Label\ru{\mu}(4,0)
\Label\lo{\nu}(0,4)
\Label\ro{\la}(4,4)
\thicklines
\Pfad(1,1),33\endPfad
\raise.5pt\hbox to 0pt{\hbox{\hskip-.5pt}\Pfad(1,1),3\endPfad\hss}%
\raise.5pt\hbox to 0pt{\hbox{\hskip-.5pt}\Pfad(2,2),3\endPfad\hss}%
\raise1pt\hbox to 0pt{\hbox{\hskip-1pt}\Pfad(1,1),3\endPfad\hss}%
\raise1pt\hbox to 0pt{\hbox{\hskip-1pt}\Pfad(2,2),3\endPfad\hss}%
\raise-.5pt\hbox to 0pt{\hbox{\hskip.5pt}\Pfad(1,1),3\endPfad\hss}%
\raise-.5pt\hbox to 0pt{\hbox{\hskip.5pt}\Pfad(2,2),3\endPfad\hss}%
\raise-1pt\hbox to 0pt{\hbox{\hskip1pt}\Pfad(1,1),3\endPfad\hss}%
\raise-1pt\hbox to 0pt{\hbox{\hskip1pt}\Pfad(2,2),3\endPfad\hss}%
\Pfad(1,3),44\endPfad
\raise.5pt\hbox to 0pt{\hbox{\hskip.5pt}\Pfad(1,3),4\endPfad\hss}%
\raise.5pt\hbox to 0pt{\hbox{\hskip.5pt}\Pfad(2,2),4\endPfad\hss}%
\raise1pt\hbox to 0pt{\hbox{\hskip1pt}\Pfad(1,3),4\endPfad\hss}%
\raise1pt\hbox to 0pt{\hbox{\hskip1pt}\Pfad(2,2),4\endPfad\hss}%
\raise-.5pt\hbox to 0pt{\hbox{\hskip-.5pt}\Pfad(1,3),4\endPfad\hss}%
\raise-.5pt\hbox to 0pt{\hbox{\hskip-.5pt}\Pfad(2,2),4\endPfad\hss}%
\raise-1pt\hbox to 0pt{\hbox{\hskip-1pt}\Pfad(1,3),4\endPfad\hss}%
\raise-1pt\hbox to 0pt{\hbox{\hskip-1pt}\Pfad(2,2),4\endPfad\hss}%
\hskip2.4cm
$$
\vskip10pt
\centerline{a. A cell without cross\hskip1.5cm
b. A cell with cross}
\caption{}
\label{fig:1}
\end{figure}

\begin{enumerate}
\item[(F1)] If $\rh=\mu=\nu$, and if there is no cross in the
cell, then $\la=\rh$.
\item[(F2)] If $\rh=\mu\ne\nu$, then $\la=\nu$.
\item[(F3)] If $\rh=\nu\ne\mu$, then $\la=\mu$.
\item[(F4)] If $\rh,\mu,\nu$ are pairwise different, 
then $\la=\mu\cup\nu$.
\item[(F5)] If $\rh\ne \mu=\nu$, 
then $\la$ is formed by adding a square to the $(k+1)$-st row of
$\mu=\nu$, given that $\mu=\nu$ and $\rh$ differ in the $k$-th row.
\item[(F6)] If $\rh=\mu=\nu$, and if there is a cross in the
cell, then $\la$ is formed by adding a square to the first row of
$\rh=\mu=\nu$.
\end{enumerate}

\begin{Remarks} (1) Due to conditions (C1) and (C2), the rules
(F1)--(F4) are forced. In particular, the uniqueness in rule (F4) is
dictated by the constraint that neighbouring partitions can only differ by
at most one square. Thus, the only ``interesting" rules are (F5)
and (F6). 

(2) Given a $0$-$1$-filling such that every row and every column
contains at most one $1$, and given a labelling of the corners along the
left side and the bottom side of a Ferrers shape by partitions such that
property (C1) and the additional property that two neighbouring
partitions along the bottom side can only be different if there is
no $1$ in the column of cells of $F$ above them, respectively two neighbouring
partitions along the left side can only be different if there is
no $1$ in the row of cells of $F$ to the right of them, the (forward) 
rules (F1)--(F6) allow one to algorithmically find the labels of
the other corners of the cells by working one's way to the right and
to the top.

(3) The rules (F5) and (F6) are carefully designed so that
one can also work one's way the other direction, that is, given
$\la,\mu,\nu$, one can reconstruct $\rh$ {\it and\/} the filling of the cell.
The corresponding ({\it backward}) {\it local rules} are:

\begin{enumerate}
\item[(B1)] If $\la=\mu=\nu$, then $\rh=\la$.
\item[(B2)] If $\la=\mu\ne\nu$, then $\rh=\nu$.
\item[(B3)] If $\la=\nu\ne\mu$, then $\rh=\mu$.
\item[(B4)] If $\la,\mu,\nu$ are pairwise different, 
then $\rh=\mu\cap\nu$.
\item[(B5)] If $\la\ne \mu=\nu$, 
then $\rh$ is formed by deleting a square from the $(k-1)$-st row of
$\mu=\nu$, given that $\mu=\nu$ and $\rh$ differ in the $k$-th row,
$k\ge2$.
\item[(B6)] If $\la\ne \mu=\nu$, and if $\la$ and $\mu=\nu$ differ in
the first row, then $\rh=\mu=\nu$.
\item[]\hskip-27pt\vbox{\hsize15.25cm\noindent
In case (B6) the cell is filled with a $1$ (an X).
In all other cases the cell is filled with a $0$.}
\end{enumerate}

Thus, given a labelling of the corners along the right/up border of a
Ferrers shape, one can algorithmically reconstruct the labels of
the other corners of the cells {\it and\/} of the $0$-$1$-filling 
by working one's way to the left and to the bottom.
\end{Remarks}

In view of the above remarks, we have the following theorem.
(See also \cite[Theorem~2.6.7]{RobyAA}, 
\cite[Theorem~3.6.3]{FomiAB}.)

\begin{Theorem} \label{thm:1}
Let $F$ be a Ferrers shape given by the $D$-$R$-sequence
$w=w_1w_2\dots w_{k}$.
The $0$-$1$-fillings of $F$ with the property that every row and
every column contains at most one $1$ are in bijection with sequences
$(\emptyset=\la^0,\la^1,\dots,\la^{k}=\emptyset)$, where $\la^{i-1}$
and $\la^i$ differ by at most one square, and
$\la^{i-1}\subseteq\la^i$ if $w_i=R$, whereas
$\la^{i-1}\supseteq\la^i$ if $w_i=D$.
Moreover, $\la^{i-1}\subsetneqq\la^i$ if and only if there is a $1$
in the column of cells of $F$ below the corners labelled by $\la^{i-1}$ and
$\la^i$, and $\la^{i-1}\supsetneqq\la^i$ if and only if there is a $1$
in the row of cells of $F$ to the left of the corners labelled by $\la^{i-1}$ and
$\la^i$.
\end{Theorem}

\begin{Remarks}(1) In analogy to classical notions 
(cf.\ \cite{RobyAA,RobyAD,SundAE}),
we call the sequences $(\la^0,\la^1,\dots,\la^{2n})$ in
this theorem {\it oscillating tableaux of type $w$
and shape $\emptyset/\emptyset$}.

(2) There are more general versions of Theorem~\ref{thm:1}
for {\it skew} Ferrers shapes. Since we have no need for these, we
refrain from reproducing them here. 
\end{Remarks}
\begin{proof}[Proof of Theorem~\ref{thm:1}]
To construct the mapping from the $0$-$1$-fillings to
the generalized oscillating tableaux, we label all the corners along
the left side and the bottom side of $F$ by $\emptyset$. Then we apply the
forward local rules (F1)--(F6) to construct labels for all the other
corners of $F$. The oscillating tableau $(\la^0,\la^1,\dots,\la^{k})$
is read off as the labels along the right/up border of $F$. For example,
the $0$-$1$-filling in Figure~\ref{fig:00}.b is mapped to the
oscillating tableau
$(\emptyset,1,1,11,11,1,1,1,\emptyset,\emptyset,1,\emptyset)$ (see
Figure~\ref{fig:0}).

That we
obtain an oscillating tableau of type $w$ is obvious from 
conditions (C1) and (C2). That the map is a bijection is obvious from
the preceding discussion.
\end{proof}

It is a now well-known fact that,
in the case that the Ferrers shape is a square and that we consider
$0$-$1$-fillings with {\it exactly} one $1$ in each row and each
column, the bijection in Theorem~\ref{thm:1} is equivalent to the
{\it Robinson--Schensted correspondence}.
Namely, $0$-$1$-fillings of an $n\times n$ square 
with exactly one $1$ in each row and each
column are in bijection with permutations. On the other hand,
according to Theorem~\ref{thm:1},
the sequences $(\la^0,\la^1,\dots,\la^{2n})$ which we read off the
top side and the right side of the square are sequences with
$\emptyset=\la^0\subsetneqq\la^1\subsetneqq\dots\subsetneqq\la^n
\supsetneqq\dots\supsetneqq\la^{2n-1}\supsetneqq\la^{2n}=\emptyset$.
In their turn, these are in bijection with pairs $(P,Q)$ of standard
tableaux of the same shape, the common shape consisting of
$n$ squares. (See \cite[Sec.~7.13]{StanBI}. The standard tableau $Q$ is defined
by the first half of the sequence, 
$\emptyset=\la^0\subsetneqq\la^1\subsetneqq\dots\subsetneqq\la^n$,
the entry $i$ being put in the square by which $\la^i$ and $\la^{i-1}$
differ, and, similarly, the standard tableau $P$ is defined
by the second half of the sequence, 
$\emptyset=\la^{2n}\subsetneqq\la^{2n-1}\subsetneqq\dots\subsetneqq\la^n$.)
It is then a theorem (cf.\ \cite[pp.~95--98]{BrFoAA},
\cite[Theorem~7.13.5]{StanBI})
that the bijection between permutations and pairs of standard tableaux
defined by the growth diagram on the square coincides with
the Robinson--Schensted correspondence, $P$ being the {\it insertion
tableau}, and $Q$ being the {\it recording tableau}.

The special case where the Ferrers shape is triangular is discussed
in more detail in the next section.

\medskip
In addition to its local description,
the bijection in Theorem~\ref{thm:1} has also a {\it global\/}
description. The latter is a consequence of a theorem of Greene
\cite{GreCAA} (see also \cite[Theorems~2.1 and 3.2]{BrFoAA}). 
In order to formulate the
result, we need the following definitions: a {\it NE-chain} of a
$0$-$1$-filling is a sequence of $1$'s in the filling such that any $1$
in the sequence is above and to the right of the preceding $1$ in 
the sequence. Similarly, a {\it SE-chain} of a
$0$-$1$-filling is a set of $1$'s in the filling such that any $1$
in the sequence is below and to the right of the preceding $1$ in the sequence.

\begin{Theorem} \label{thm:3}
Given a growth diagram with empty partitions labelling all the
corners along the left side and the bottom side of the Ferrers shape, 
the partition $\la=(\la_1,\la_2,\break\dots,\la_\ell)$ labelling corner $c$ 
satisfies the following two properties:
\begin{enumerate}
\item [(G1)]For any $k$, the maximal cardinality of the union of $k$
NE-chains situated in the rectangular region to the left and
below of $c$ is equal to $\la_1+\la_2+\dots+\la_k$.
\item [(G2)]For any $k$, the maximal cardinality of the union of $k$
SE-chains situated in the rectangular region to the left and
below of $c$ is equal to $\la'_1+\la'_2+\dots+\la'_k$, where $\la'$
denotes the partition conjugate to $\la$.
\end{enumerate}
In particular, $\la_1$ is the length of the longest NE-chain
in the rectangular region to the left and below of $c$, and $\la'_1$
is the length of the longest SE-chain in the same
rectangular region.
\end{Theorem}

In order to formulate the main theorem of this section, we introduce
one more piece of notation. We write $N(F;n;\NE=s,\SE=t)$ for the number
of $0$-$1$-fillings of the Ferrers shape $F$ with exactly $n$ $1$'s,
such that there is at most one $1$ in each column and in each row,
and such that the longest
NE-chain has length $s$ and the longest SE-chain, the
smallest rectangle containing the chain being contained in $F$, has
length $t$.

\begin{Theorem} \label{thm:2}
For any Ferrers shape $F$ and positive integers $s$ and $t$, we have
$$N(F;n;\NE=s,\SE=t)=N(F;n;\NE=t,\SE=s).$$
\end{Theorem}
\begin{proof}We define a bijection between the $0$-$1$-fillings
counted by $N(F;n;\NE=s,\SE=t)$ and those counted by $N(F;n;\NE=t,\SE=s)$.
Let the Ferrers shape $F$ be given by the $D$-$R$-sequence
$w=w_1w_2\dots w_k$.
Given a $0$-$1$-filling counted by $N(F;n;\NE=s,\SE=t)$ we apply the
mapping of the proof of Theorem~\ref{thm:1}. Thus, we obtain an
oscillating tableau
$(\emptyset=\la^0,\la^1,\dots,\la^{k}=\emptyset)$. 
Since the $0$-$1$-filling had $n$ entries $1$, in the oscillating
tableau there are exactly $n$ ``rises" $\la^{i-1}\subsetneqq\la^i$
(and, hence, exactly $n$ ``falls" $\la^{i-1}\supsetneqq\la^i$).
Moreover, by Theorem~\ref{thm:3}, we have $\la^i_1\le s$ and
$(\la^j)'_1\le t$ for all $i$ and $j$, with equality for at least one
$i$ and at least one $j$. Now we apply the inverse mapping to
the sequence $(\emptyset=(\la^0)',(\la^1)',\dots,(\la^{k})'=\emptyset)$
of conjugate partitions. Thus, we obtain a $0$-$1$-filling counted
by $N(F;n;\NE=t,\SE=s)$.
\end{proof}

If we specialize Theorem~\ref{thm:2} to the case where $F$ is a
square and to $0$-$1$-fillings which have exactly one entry $1$ in
each row and in each column, then we obtain a trivial statement:
The number of permutations of $\{1,2,\dots,n\}$ with longest increasing
subsequence of length $s$ and longest decreasing
subsequence of length $t$ is equal to the 
number of permutations of $\{1,2,\dots,n\}$ with longest increasing
subsequence of length $t$ and longest decreasing
subsequence of length $s$. 
This statement is indeed trivial because, given a former permutation 
$\pi_1\pi_2\dots\pi_n$, the reversal $\pi_n\pi_{n-1}\dots\pi_1$
will belong to the latter permutations.

However, for other Ferrers shapes $F$, there is no trivial
explanation for Theorem~\ref{thm:2}. In particular,
as we are going to show in the next section,
if we specialize Theorem~\ref{thm:2} to the case where $F$ is 
triangular, then we obtain
a non-obvious theorem, originally due to Chen et al.\ \cite{ChDDAB}, 
about crossings and nestings in set partitions.

\section{Crossings and nestings in set partitions and matchings}
\label{sec:part}

In this section we show that the main theorems in \cite{ChDDAB} are,
essentially, special cases of the theorems in the previous section.

First of all, we have to recall the definitions and basic objects
from \cite{ChDDAB}. The objects of consideration in \cite{ChDDAB} are
{\it set partitions} of $\{1,2,\dots,n\}$. A block $\{i_1,i_2,\dots,i_b\}$,
$i_1<i_2<\dots<i_b$, of such a set partition is represented by
the set of pairs $\{(i_1,i_2),(i_2,i_3),\dots,\break(i_{b-1},i_b)\}$.
More generally, a set partition is represented by the union 
of all sets of pairs, the union being taken over all its blocks. 
This representation is called the {\it
standard representation} of the set partition.
For example, the set partition
$\{\{1,4,5,7\},\{2,6\},\{3\}\}$ is represented as the set
$\{(1,4),(4,5),(5,7),(2,6)\}$. Next, one defines a {\it $k$-crossing} of a
set partition to be a subset $\{(i_1,j_1),(i_2,j_2),\dots,(i_k,j_k)\}$ of
its standard representation where $i_1<i_2<\dots<i_k<j_1<j_2<\dots<j_k$. 
Similarly, one defines a {\it $k$-nesting} of a set partition 
to be a subset $\{(i_1,j_1),(i_2,j_2),\dots,(i_k,j_k)\}$ of
its standard representation where
$i_1<i_2<\dots<i_k<j_k<\dots<j_2<j_1$. (These notions have 
intuitive pictorial meanings if one connects a pair $(i,j)$ in the
standard representation of a set partition by an arc, cf.\ \cite{ChDDAB}.)
Finally, given a set partition $P$, we write $\cross(P)$ for the maximal
number $k$ such that $P$ has a $k$-crossing, and we write $\nest(P)$ 
for the maximal number $k$ such that $P$ has a $k$-nesting,

In \cite[Eq.~(3)]{ChDDAB}, the following theorem is proved
by a Robinson--Schensted-like insertion/deletion process, which sets
up a bijection between set partitions and ``vacillating tableaux"
(see below).

\begin{Theorem} \label{thm:4}
Let $n,s,t$ be positive integers. Then
the number of set partitions of $\{1,2,\dots,n\}$ with $\cross(P)=s$ and
$\nest(P)=t$ is equal to the 
number of set partitions of $\{1,2,\dots,n\}$ with $\cross(P)=t$ and
$\nest(P)=s$.
\end{Theorem}

\begin{figure}[h]
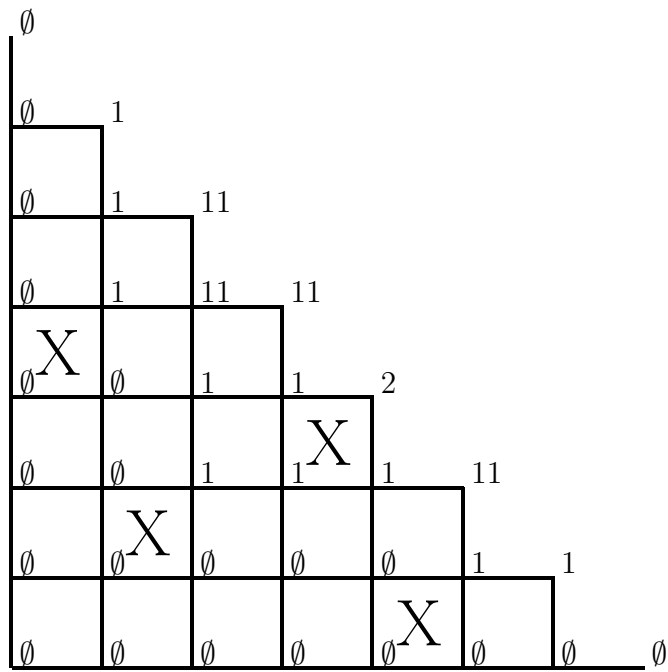

$$
\Einheit.4cm
\Pfad(0,21),666666666666666666666111111111111111111111\endPfad
\Pfad(0,18),111666666666666666666\endPfad
\Pfad(0,15),111111666666666666666\endPfad
\Pfad(0,12),111111111666666666666\endPfad
\Pfad(0,9),111111111111666666666\endPfad
\Pfad(0,6),111111111111111666666\endPfad
\Pfad(0,3),111111111111111111666\endPfad
\Label\ro{\emptyset}(0,0)
\Label\ro{\emptyset}(0,3)
\Label\ro{\emptyset}(0,6)
\Label\ro{\emptyset}(0,9)
\Label\ro{\emptyset}(0,12)
\Label\ro{\emptyset}(0,15)
\Label\ro{\emptyset}(0,18)
\Label\ro{\emptyset}(0,21)
\Label\ro{1}(3,18)
\Label\ro{1}(3,15)
\Label\ro{1}(3,12)
\Label\ro{\emptyset}(3,9)
\Label\ro{\emptyset}(3,6)
\Label\ro{\emptyset}(3,3)
\Label\ro{\emptyset}(3,0)
\Label\ro{\hphantom{1}11}(6,15)
\Label\ro{\hphantom{1}11}(6,12)
\Label\ro{1}(6,9)
\Label\ro{1}(6,6)
\Label\ro{\emptyset}(6,3)
\Label\ro{\emptyset}(6,0)
\Label\ro{\hphantom{1}11}(9,12)
\Label\ro{1}(9,9)
\Label\ro{1}(9,6)
\Label\ro{\emptyset}(9,3)
\Label\ro{\emptyset}(9,0)
\Label\ro{2}(12,9)
\Label\ro{1}(12,6)
\Label\ro{\emptyset}(12,3)
\Label\ro{\emptyset}(12,0)
\Label\ro{\hphantom{2}11}(15,6)
\Label\ro{1}(15,3)
\Label\ro{\emptyset}(15,0)
\Label\ro{1}(18,3)
\Label\ro{\emptyset}(18,0)
\Label\ro{\emptyset}(21,0)
\Label\ro{\text {\Huge X}}(4,4)
\Label\ro{\text {\Huge X}}(1,10)
\Label\ro{\text {\Huge X}}(13,1)
\Label\ro{\text {\Huge X}}(10,7)
\hskip8.4cm
$$
\caption{Bijection between $\{\{1,4,5,7\},\{2,6\},\{3\}\}$
and\newline
$(\emptyset,\emptyset,1,1,11,11,11,1,2,1,11,1,1,\emptyset,\emptyset)$}
\label{fig:3}
\end{figure}

As we now explain, this is just a special case of
Theorem~\ref{thm:2}, where $F$ is triangular.
More precisely, let $\Delta_n$ be the triangular shape with $n-1$ cells in
the bottom row, $n-2$ cells in the row above, etc., and $1$ cell in
the top-most row. See Figure~\ref{fig:3} for an example in which
$n=7$. (The filling and labelling of the corners should be ignored
at this point. For convenience, we also joined pending edges at the
right and at the top of $\Delta_n$.)

We represent a set partition of $\{1,2,\dots,n\}$, given by
its standard representation
$\{(i_1,j_1),(i_2,j_2),\dots,(i_m,j_m)\}$,
as a $0$-$1$-filling, by putting a $1$ in the $i_r$-th column and the
$j_r$-th row from above (where we number rows such that the
row consisting of $j-1$ cells is the row numbered $j$), $r=1,2,\dots,m$. 
The filling corresponding to the set partition
$\{\{1,4,5,7\},\{2,6\},\{3\}\}$ is shown in Figure~\ref{fig:3},
where we present again $1$'s by X's and suppress the $0$'s.

It is obvious that this defines a bijection between set partitions
of $\{1,2,\dots,n\}$ and $0$-$1$-fillings of $\Delta_n$ in which every
row and every column contains at most one $1$. Moreover, a
$k$-crossing corresponds to a SE-chain in the fillings, while a
$k$-nesting corresponds to a NE-chain. Thus, Theorem~\ref{thm:2}
specialized to $F=\Delta_n$ yields Theorem~\ref{thm:4} immediately.

Figure~\ref{fig:3} shows as well the labelling of the corners
by partitions if we apply the correspondence of Theorem~\ref{thm:1}.
The ``extra" corners created by the pending edges are labelled
by empty partitions.
The oscillating tableaux which one reads along the right/up border of
$\Delta_n$ are sequences
$(\emptyset=\la^0,\la^1,\dots,\la^{2n}=\emptyset)$ with the property
that $\la^{2i+1}$ is obtained from $\la^{2i}$ by doing nothing or by
deleting one square, and $\la^{2i}$ is obtained from $\la^{2i-1}$ by
doing nothing or adding a square. Such oscillating tableaux are
called {\it vacillating tableaux} in \cite{ChDDAB}. Again, from
Theorem~\ref{thm:1}, it is obvious that set partitions of
$\{1,2,\dots,n\}$ are in bijection with vacillating tableaux 
$(\emptyset=\la^0,\la^1,\dots,\la^{2n}=\emptyset)$.
Figure~\ref{fig:3} does in fact work out the bijection for
the set partition of Example~4 in \cite{ChDDAB}. One can verify
that the resulting vacillating tableau is the same as in
\cite{ChDDAB}. Indeed, in general, the bijection between set partitions
and vacillating tableaux described in Section~3 of \cite{ChDDAB}
is equivalent with our growth diagram bijection. This follows
from the afore-mentioned 
fact that the growth diagrams model Robinson--Schensted
insertion.

The main theorem in \cite{ChDDAB} is in fact a refinement of
Theorem~\ref{thm:4}. This refinement takes also into account the
minimal and the maximal elements in the blocks of a set partition $P$.
While, from the growth diagram point of view, this refinement
is special for triangular Ferrers shapes and cannot be extended in a
natural way to arbitrary Ferrers shapes, it follows nevertheless with
equal ease from the growth diagram point of view. 
Given a set partition $P$, let $\min(P)$ be the set of minimal
elements of the blocks of $P$, and let $\max(P)$ be the set of maximal
elements of the blocks of $P$. Then Theorem~1 from \cite{ChDDAB}
reads as follows.

\begin{Theorem} \label{thm:5}
Let $n,s,t$ be positive integers, and let $S$ and $T$ be two subsets
of $\{1,2,\break\dots,n\}$. Then
the number of set partitions of $\{1,2,\dots,n\}$ with $\cross(P)=s$,
$\nest(P)=t$, $\min(P)=S$, $\max(P)=T$ is equal to the 
number of set partitions of $\{1,2,\dots,n\}$ with $\cross(P)=t$,
$\nest(P)=s$, $\min(P)=S$, $\max(P)=T$.
\end{Theorem}
\begin{proof} As in the proof of Theorem~\ref{thm:2}, we set up
a bijection between the two different sets of set partitions.
In fact, using the correspondence between set
partitions and $0$-$1$-fillings that we explained after the statement
of Theorem~\ref{thm:4}, the bijection is exactly the same as the one
in the proof of Theorem~\ref{thm:2}. We only have to figure out
how one can detect the minimal and maximal elements in blocks
in the vacillating tableau to which a set partition is mapped,
and verify that these are kept invariant under conjugation of
the partitions of the vacillating tableau. 

Let $\{i_1,i_2,\dots,i_b\}$ be a block of the set partition $P$.
Then, in the $0$-$1$-filling corresponding to $P$, there is a $1$ in
{\it column} $i_1$ and row $i_2$, while there is no $1$ in {\it row} $i_1$.
Similarly, there is a $1$ in
column $i_{b-1}$ and {\it row} $i_b$, while there is no $1$ in {\it column} $i_b$.
Consequently, the two corners on the right of row $i_1$ will be labelled
by partitions $\la^{2i_1-2}$ and $\la^{2i_1-1}$ with
$\la^{2i_1-2}=\la^{2i_1-1}$, and 
the two corners on the top of column $i_1$ will be labelled
by partitions $\la^{2i_1-1}$ and $\la^{2i_1}$ with
$\la^{2i_1-1}\subsetneqq\la^{2i_1}$. On the other hand,
the two corners on the right of row $i_b$ will be labelled
by partitions $\la^{2i_b-2}$ and $\la^{2i_b-1}$ with
$\la^{2i_b-2}\supsetneqq\la^{2i_b-1}$, and 
the two corners on the top of column $i_b$ will be labelled
by partitions $\la^{2i_b-1}$ and $\la^{2i_b}$ with
$\la^{2i_b-1}=\la^{2i_b}$. Thus, in summary, the sets $\min(P)$ and
$\max(P)$ can be detected from the growth properties of
the subsequences $(\la^{2i-2},\la^{2i-1},\la^{2i})$ of the
vacillating tableau corresponding to $P$. Clearly, these remain
invariant under conjugation of the partitions. Hence, the theorem.
\end{proof}

For the sake of completeness, we also review how one can realize
the other Robinson--Schensted-like insertion/deletion mappings 
in \cite{ChDDAB}. 

\begin{figure}[h]
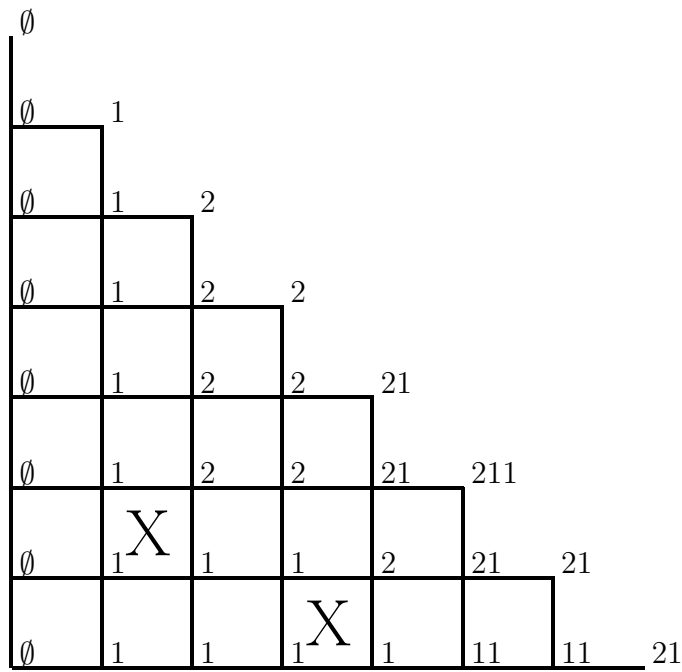

$$
\Einheit.4cm
\Pfad(0,21),666666666666666666666111111111111111111111\endPfad
\Pfad(0,18),111666666666666666666\endPfad
\Pfad(0,15),111111666666666666666\endPfad
\Pfad(0,12),111111111666666666666\endPfad
\Pfad(0,9),111111111111666666666\endPfad
\Pfad(0,6),111111111111111666666\endPfad
\Pfad(0,3),111111111111111111666\endPfad
\Label\ro{\emptyset}(0,0)
\Label\ro{\emptyset}(0,3)
\Label\ro{\emptyset}(0,6)
\Label\ro{\emptyset}(0,9)
\Label\ro{\emptyset}(0,12)
\Label\ro{\emptyset}(0,15)
\Label\ro{\emptyset}(0,18)
\Label\ro{\emptyset}(0,21)
\Label\ro{1}(3,18)
\Label\ro{1}(3,15)
\Label\ro{1}(3,12)
\Label\ro{1}(3,9)
\Label\ro{1}(3,6)
\Label\ro{1}(3,3)
\Label\ro{1}(3,0)
\Label\ro{2}(6,15)
\Label\ro{2}(6,12)
\Label\ro{2}(6,9)
\Label\ro{2}(6,6)
\Label\ro{1}(6,3)
\Label\ro{1}(6,0)
\Label\ro{2}(9,12)
\Label\ro{2}(9,9)
\Label\ro{2}(9,6)
\Label\ro{1}(9,3)
\Label\ro{1}(9,0)
\Label\ro{\hphantom{2}21}(12,9)
\Label\ro{\hphantom{2}21}(12,6)
\Label\ro{2}(12,3)
\Label\ro{1}(12,0)
\Label\ro{\hphantom{22}211}(15,6)
\Label\ro{\hphantom{2}21}(15,3)
\Label\ro{\hphantom{2}11}(15,0)
\Label\ro{\hphantom{2}21}(18,3)
\Label\ro{\hphantom{2}11}(18,0)
\Label\ro{\hphantom{2}21}(21,0)
\Label\ro{\text {\Huge X}}(4,4)
\Label\ro{\text {\Huge X}}(10,1)
\hskip8.4cm
$$
\caption{Bijection between the pair in \eqref{eq:PT}
and the vacillating tableau in \eqref{eq:vac}}
\label{fig:2}
\end{figure}

To begin with, Theorem~5 in \cite{ChDDAB} describes a bijection
between pairs $(P,T)$ of set partitions $P$ of $\{1,2,\dots,n\}$ 
and (partial) standard 
tableaux $T$ of shape $\la$ with set of entries contained in $\max(P)$ and
vacillating tableaux $(\emptyset=\la^0,\la^1,\dots,\la^{2n}=\la)$.
In growth diagram language, the bijection can be realized by
putting the $0$-$1$-filling corresponding to $P$ in the cells of
$\Delta_n$, labelling all the corners along the left side of $\Delta_n$
by $\emptyset$, and labelling the corners along the bottom side of
$\Delta_n$ by the increasing sequence of partitions corresponding to
$T$, that is, the $\ell$-th corner is labelled by the partition
corresponding to the shape that is covered by the entries of $T$ contained
in $\{1,2,\dots,\ell\}$, $\ell=0,1,\dots,n$. Figure~\ref{fig:2} shows 
the growth diagram describing the bijection in the case that 
\begin{equation} \label{eq:PT}
(P,T)=\left(\begin{matrix} 1\ 7\\5\
\hphantom{3}\end{matrix}\ ,\
\{\{1\},\{2,6\},\{3\},\{4,7\},\{5\}\}\right).
\end{equation}
This is, in fact, the pair $(P,T)$ of Example~3 in \cite{ChDDAB}.
As the figure shows, this pair is mapped to the vacillating tableau
\begin{equation} \label{eq:vac}
(\emptyset,\emptyset,1,1,2,2,2,2,21,21,211,21,21,11,21),
\end{equation}
in agreement with the result of the (differently defined) bijection
in \cite{ChDDAB}.

In Section~4 of \cite{ChDDAB}, a variant of the bijection between
set partitions and vacillating tableaux is discussed, namely
a bijection between set partitions and {\it hesitating tableaux}.
Here, a hesitating tableau is an oscillating tableau
$(\emptyset=\la^0,\la^1,\dots,\la^{2n}=\emptyset)$
with the property
that, for each $i$, either 
(1) $\la^{2i-2}=\la^{2i-1}\subsetneqq\la^{2i}$, or
(2) $\la^{2i-2}\supsetneqq\la^{2i-1}=\la^{2i}$, or
(3) $\la^{2i-1}\subsetneqq\la^{2i-2}\supsetneqq\la^{2i}$.

\begin{figure}[h]
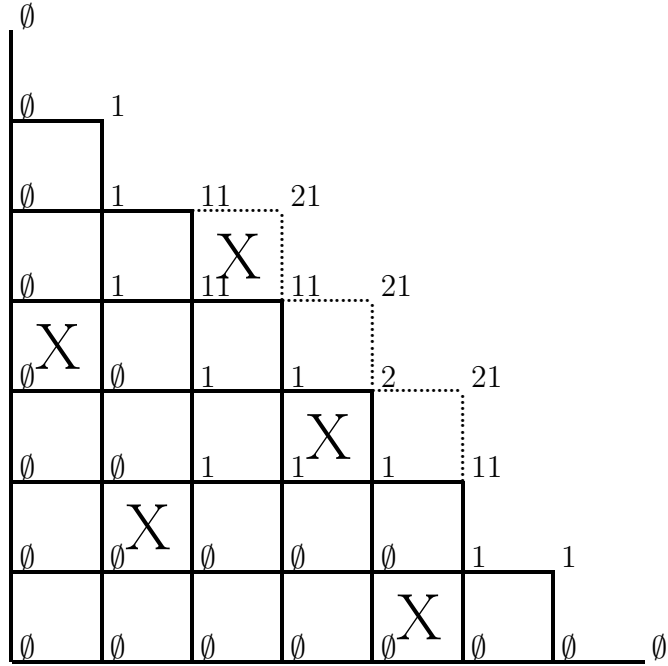

$$
\Einheit.4cm
\Pfad(0,21),666666666666666666666111111111111111111111\endPfad
\SPfad(6,15),111666111666111666\endSPfad
\Pfad(0,18),111666666666666666666\endPfad
\Pfad(0,15),111111666666666666666\endPfad
\Pfad(0,12),111111111666666666666\endPfad
\Pfad(0,9),111111111111666666666\endPfad
\Pfad(0,6),111111111111111666666\endPfad
\Pfad(0,3),111111111111111111666\endPfad
\Label\ro{\emptyset}(0,0)
\Label\ro{\emptyset}(0,3)
\Label\ro{\emptyset}(0,6)
\Label\ro{\emptyset}(0,9)
\Label\ro{\emptyset}(0,12)
\Label\ro{\emptyset}(0,15)
\Label\ro{\emptyset}(0,18)
\Label\ro{\emptyset}(0,21)
\Label\ro{1}(3,18)
\Label\ro{1}(3,15)
\Label\ro{1}(3,12)
\Label\ro{\emptyset}(3,9)
\Label\ro{\emptyset}(3,6)
\Label\ro{\emptyset}(3,3)
\Label\ro{\emptyset}(3,0)
\Label\ro{\hphantom{1}11}(6,15)
\Label\ro{\hphantom{1}11}(6,12)
\Label\ro{1}(6,9)
\Label\ro{1}(6,6)
\Label\ro{\emptyset}(6,3)
\Label\ro{\emptyset}(6,0)
\Label\ro{\hphantom{1}11}(9,12)
\Label\ro{1}(9,9)
\Label\ro{1}(9,6)
\Label\ro{\emptyset}(9,3)
\Label\ro{\emptyset}(9,0)
\Label\ro{2}(12,9)
\Label\ro{1}(12,6)
\Label\ro{\emptyset}(12,3)
\Label\ro{\emptyset}(12,0)
\Label\ro{\hphantom{2}11}(15,6)
\Label\ro{1}(15,3)
\Label\ro{\emptyset}(15,0)
\Label\ro{1}(18,3)
\Label\ro{\emptyset}(18,0)
\Label\ro{\emptyset}(21,0)
\Label\ro{\hphantom{1}21}(9,15)
\Label\ro{\hphantom{1}21}(12,12)
\Label\ro{\hphantom{1}21}(15,9)
\Label\ro{\text {\Huge X}}(4,4)
\Label\ro{\text {\Huge X}}(1,10)
\Label\ro{\text {\Huge X}}(13,1)
\Label\ro{\text {\Huge X}}(10,7)
\Label\ro{\text {\Huge X}}(7,13)
\hskip8.4cm
$$
\caption{Bijection between $\{\{1,4,5,7\},\{2,6\},\{3\}\}$
and\newline
$(\emptyset,\emptyset,1,1,11,21,11,21,2,21,11,1,1,\emptyset,\emptyset)$}
\label{fig:4}
\end{figure}

In terms of growth diagrams, this bijection can be again realized as
a special case of the bijection described in the proof of
Theorem~\ref{thm:2}. (This realization is also attributed to
Michael Korn in \cite{ChDDAB}.) Here, we do not transform set
partitions into $0$-$1$-fillings of a triangular shape, but 
into $0$-$1$-fillings of a slightly modified shape
that may vary depending on the original set partition. 
More precisely, 
given a set partition $P$ of $\{1,2,\dots,n\}$, we deform
$\Delta_n$ by attaching an extra cell in the $i$-th row from above and 
the $i$-th column for every singleton block $\{i\}$ of $P$,
and by attaching an extra cell in the $j$-th row from above and 
the $j$-th column whenever $(i,j)$ and $(j,k)$ are both in the
standard representation of $P$, for some $i$ and $k$.
Then, as before, we transform $P$
into a $0$-$1$-filling of this modification of $\Delta_n$, 
by placing a $1$ in the $i$-th column and $j$-th row if $(i,j)$
is a pair in the standard representation of $P$, and, in addition,
by placing a $1$ in the added cell in the $i$-th row and 
the $i$-th column for every singleton block $\{i\}$ of $P$.
All other cells are filled with $0$'s.
An example is shown in Figure~\ref{fig:4}. It shows the
$0$-$1$-filling corresponding to the partition
$\{\{1,4,5,7\},\{2,6\},\{3\}\}$ under this modified rule.
The added cells are indicated by dotted lines.

To realize the bijection to hesitating tableaux, we label all the
corners along the left side and the bottom side by $\emptyset$,
and then apply the forward local rules (F1)--(F6) to determine the
labels for all other corners. Along the right/up border one reads off
a hesitating tableau. It is easy to see that this defines a
bijection. The example in Figure~\ref{fig:4} maps the partition
$\{\{1,4,5,7\},\{2,6\},\{3\}\}$ to the hesitating tableau
$(\emptyset,\emptyset,1,1,11,21,11,21,2,21,11,1,1,\emptyset,\emptyset)$.
This example corresponds to Example~6 in \cite{ChDDAB}.

For the sake of completeness, we record the consequence of this
bijection from \cite[Theorem~11]{ChDDAB}. In the statement,
we need the notion of {\it enhanced\/} $k$-crossings and
$k$-nestings. To define these, one first defines the {\it enhanced
representation} of a set partition $P$ to be union of the standard
representation of $P$ with the set of pairs $(i,i)$, where $i$ ranges
over all the singleton blocks $\{i\}$ of $P$. Then
one defines an {\it enhanced $k$-crossing} of a
set partition to be a subset $\{(i_1,j_1),(i_2,j_2),\dots,(i_k,j_k)\}$ of
its enhanced representation where $i_1<i_2<\dots<i_k\le j_1<j_2<\dots<j_k$. 
Similarly, one defines an {\it enhanced $k$-nesting} of a set partition 
to be a subset $\{(i_1,j_1),(i_2,j_2),\dots,(i_k,j_k)\}$ of
its enhanced representation where
$i_1<i_2<\dots<i_k\le j_k<\dots<j_2<j_1$.
Finally, given a set partition $P$, we write $\ccross(P)$ for the maximal
number $k$ such that $P$ has an enhanced $k$-crossing, and we write
$\nnest(P)$ for the maximal number $k$ such that $P$ has an
enhanced $k$-nesting,

With the above notation, the following variant of Theorem~\ref{thm:5}
holds.

\begin{Theorem} \label{thm:6}
Let $n,s,t$ be positive integers, and let $S$ and $T$ be two subsets
of $\{1,2,\break\dots,n\}$. Then
the number of set partitions of $\{1,2,\dots,n\}$ with $\ccross(P)=s$,
$\nnest(P)=t$, $\min(P)=S$, $\max(P)=T$ is equal to the 
number of set partitions of $\{1,2,\dots,n\}$ with $\ccross(P)=t$,
$\nnest(P)=s$, $\min(P)=S$, $\max(P)=T$.
\end{Theorem}

\begin{figure}[h]
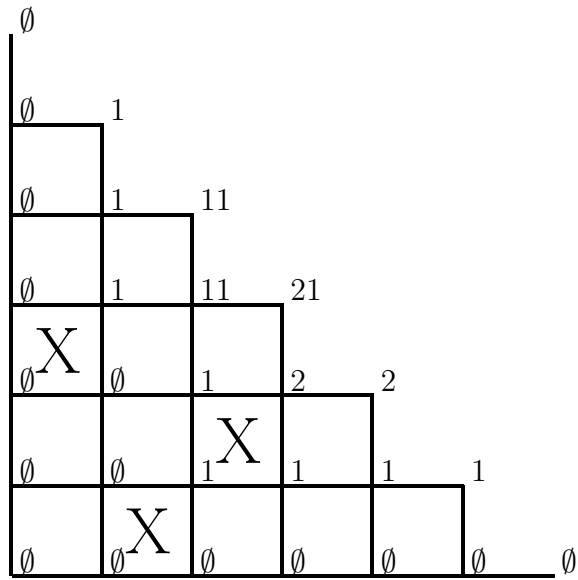

$$
\Einheit.4cm
\Pfad(0,18),666666666666666666111111111111111111\endPfad
\Pfad(0,15),111666666666666666\endPfad
\Pfad(0,12),111111666666666666\endPfad
\Pfad(0,9),111111111666666666\endPfad
\Pfad(0,6),111111111111666666\endPfad
\Pfad(0,3),111111111111111666\endPfad
\Label\ro{\emptyset}(0,0)
\Label\ro{\emptyset}(0,3)
\Label\ro{\emptyset}(0,6)
\Label\ro{\emptyset}(0,9)
\Label\ro{\emptyset}(0,12)
\Label\ro{\emptyset}(0,15)
\Label\ro{\emptyset}(0,18)
\Label\ro{1}(3,15)
\Label\ro{1}(3,12)
\Label\ro{1}(3,9)
\Label\ro{\emptyset}(3,6)
\Label\ro{\emptyset}(3,3)
\Label\ro{\emptyset}(3,0)
\Label\ro{\hphantom{1}11}(6,12)
\Label\ro{\hphantom{1}11}(6,9)
\Label\ro{1}(6,6)
\Label\ro{1}(6,3)
\Label\ro{\emptyset}(6,0)
\Label\ro{\hphantom{1}21}(9,9)
\Label\ro{2}(9,6)
\Label\ro{1}(9,3)
\Label\ro{\emptyset}(9,0)
\Label\ro{2}(12,6)
\Label\ro{1}(12,3)
\Label\ro{\emptyset}(12,0)
\Label\ro{1}(15,3)
\Label\ro{\emptyset}(15,0)
\Label\ro{\emptyset}(18,0)
\Label\ro{\text {\Huge X}}(4,1)
\Label\ro{\text {\Huge X}}(1,7)
\Label\ro{\text {\Huge X}}(7,4)
\hskip8.4cm
$$
\caption{Bijection between $\{\{1,4\},\{2,6\},\{3,5\}\}$ and
$\{\emptyset,1,11,21,2,1,\emptyset\}$}
\label{fig:5}
\end{figure}

We conclude this section by recalling the growth diagram bijection
between {\it matchings} and (ordinary) oscillating tableaux
(cf.\ \cite[Sec.~4]{RobyAA} or \cite{RobyAD}).
Clearly, matchings of $\{1,2,\dots,2n\}$ can be alternatively seen as
partitions of $\{1,2,\dots,2n\}$ all the blocks of which consist of
two elements. In their turn, if we transform the latter to the
corresponding $0$-$1$-fillings of the triangular shape 
$\Delta_n$, then these filling have the property that in the union of
the $i$-th column and the $i$-th row (from above) there is exactly
one $1$. See Figure~\ref{fig:5} for the $0$-$1$-filling corresponding
to the matching $\{\{1,4\},\{2,6\},\{3,5\}\}$. Consequently, if we
apply the bijection from the proof of Theorem~\ref{thm:2} we obtain a
vacillating tableau 
$(\emptyset=\la^0,\la^1,\dots,\la^{4n}=\emptyset)$
with the property
that, for each $i$, either 
(1) $\la^{2i-2}=\la^{2i-1}\subsetneqq\la^{2i}$, or
(2) $\la^{2i-2}\supsetneqq\la^{2i-1}=\la^{2i}$. 
Thus, the information
contained in the partitions $\la^{2i-1}$ is superfluous and can be
dropped. What remains is an (ordinary) oscillating tableau,
that is, a sequence
$(\emptyset=\la^0,\la^2,\la^4\dots,\la^{4n}=\emptyset)$
with the property that successive partitions in the sequence differ
by exactly one square. Figure~\ref{fig:5} shows that, under this bijection,
the matching $\{\{1,4\},\{2,6\},\{3,5\}\}$ is mapped to the
oscillating tableau $\{\emptyset,1,11,21,2,1,\emptyset\}$.
This is in accordance with Example~7 in \cite{ChDDAB}.

\section{Growth diagrams for arbitrary fillings of Ferrers shapes}\label{sec:RSK}

In this section we embed the considerations of the previous two
sections into the larger context where we relax the conditions on the
fillings that we imposed so far: namely, in this section, 
we do not insist anymore that our fillings have at most one $1$ in
each row and in each column and otherwise $0$'s. That is, we shall
allow more than one $1$ in rows and columns and we shall also allow
arbitrary non-negative entries in our fillings.
As it turns out,
there are now four variants how to define growth diagrams for these
more general fillings which lead to (different) extensions of
Theorem~\ref{thm:2}. (All of these are, of course, special instances
of the general set-up in \cite[Theorem~3.6]{FomiAF}.) 
They have been described in some
detail in \cite{RobyAA} and \cite{LeeuAH}, albeit without giving the
analogues of Greene's theorem in all four cases, the latter being
exactly what we need for our purposes.
We use the opportunity here to provide thorough descriptions of all
four of these correspondences in this section, including the
analogues of Greene's theorem. 
The corresponding extensions of Theorems~\ref{thm:1} and \ref{thm:3} 
are given in Theorems~\ref{thm:1a}--\ref{thm:3d},
while the corresponding extensions of Theorem~\ref{thm:2} are given
in Theorem~\ref{thm:2a}.

\begin{figure}[h]
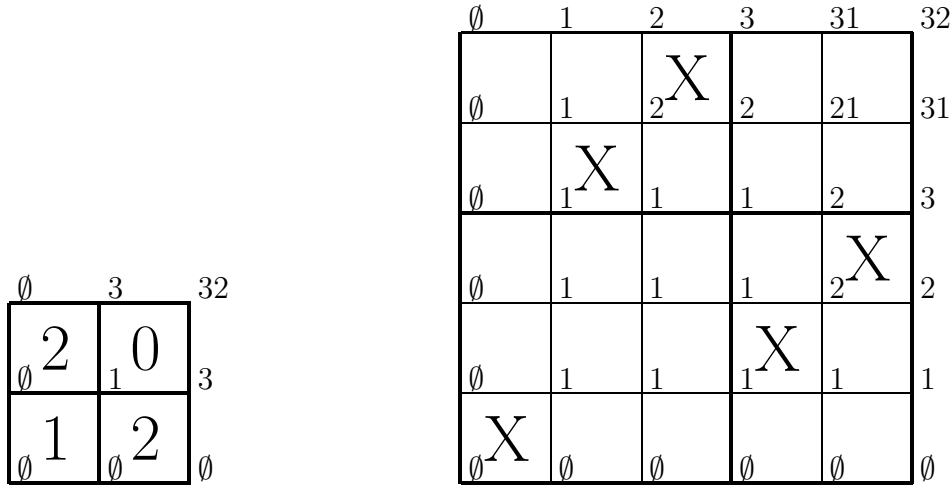

$$
\Einheit.4cm
\Pfad(0,6),111111\endPfad
\Pfad(0,3),111111\endPfad
\Pfad(0,0),111111\endPfad
\Pfad(0,0),222222\endPfad
\Pfad(3,0),222222\endPfad
\Pfad(6,0),222222\endPfad
\Label\ro{\text {\Huge $1$}}(1,1)
\Label\ro{\text {\Huge $2$}}(1,4)
\Label\ro{\text {\Huge $2$}}(4,1)
\Label\ro{\text {\Huge $0$}}(4,4)
\Label\ro{\emptyset}(0,0)
\Label\ro{\emptyset}(0,3)
\Label\ro{\emptyset}(0,6)
\Label\ro{\emptyset}(3,0)
\Label\ro{1}(3,3)
\Label\ro{3}(3,6)
\Label\ro{\emptyset}(6,0)
\Label\ro{3}(6,3)
\Label\ro{\hphantom{1}32}(6,6)
\hbox{\hskip6cm}
\PfadDicke{1pt}
\Pfad(0,15),111111111111111\endPfad
\PfadDicke{.3pt}
\Pfad(0,12),111111111111111\endPfad
\PfadDicke{1pt}
\Pfad(0,9),111111111111111\endPfad
\PfadDicke{.3pt}
\Pfad(0,6),111111111111111\endPfad
\Pfad(0,3),111111111111111\endPfad
\PfadDicke{1pt}
\Pfad(0,0),111111111111111\endPfad
\Pfad(15,0),222222222222222\endPfad
\PfadDicke{.3pt}
\Pfad(12,0),222222222222222\endPfad
\PfadDicke{1pt}
\Pfad(9,0),222222222222222\endPfad
\PfadDicke{.3pt}
\Pfad(6,0),222222222222222\endPfad
\Pfad(3,0),222222222222222\endPfad
\PfadDicke{1pt}
\Pfad(0,0),222222222222222\endPfad
\Label\ro{\emptyset}(0,0)
\Label\ro{\emptyset}(0,3)
\Label\ro{\emptyset}(0,6)
\Label\ro{\emptyset}(0,9)
\Label\ro{\emptyset}(0,12)
\Label\ro{\emptyset}(0,15)
\Label\ro{1}(3,15)
\Label\ro{1}(3,12)
\Label\ro{1}(3,9)
\Label\ro{1}(3,6)
\Label\ro{1}(3,3)
\Label\ro{\emptyset}(3,0)
\Label\ro{2}(6,15)
\Label\ro{2}(6,12)
\Label\ro{1}(6,9)
\Label\ro{1}(6,6)
\Label\ro{1}(6,3)
\Label\ro{\emptyset}(6,0)
\Label\ro{3}(9,15)
\Label\ro{2}(9,12)
\Label\ro{1}(9,9)
\Label\ro{1}(9,6)
\Label\ro{1}(9,3)
\Label\ro{\emptyset}(9,0)
\Label\ro{\hphantom{1}31}(12,15)
\Label\ro{\hphantom{1}21}(12,12)
\Label\ro{2}(12,9)
\Label\ro{2}(12,6)
\Label\ro{1}(12,3)
\Label\ro{\emptyset}(12,0)
\Label\ro{\hphantom{1}32}(15,15)
\Label\ro{\hphantom{1}31}(15,12)
\Label\ro{3}(15,9)
\Label\ro{2}(15,6)
\Label\ro{1}(15,3)
\Label\ro{\emptyset}(15,0)
\Label\ro{\text {\Huge X}}(1,1)
\Label\ro{\text {\Huge X}}(4,10)
\Label\ro{\text {\Huge X}}(7,13)
\Label\ro{\text {\Huge X}}(10,4)
\Label\ro{\text {\Huge X}}(13,7)
\hskip6cm
$$
\caption{The ``blow-up" of an arbitrary filling}
\label{fig:6a}
\end{figure}

To sketch the idea, let us consider the filling of the $2\times 2$
square on the left of Figure~\ref{fig:6a}. 
(At this point, the labellings of the
corners of the cells should be ignored.)
We cannot apply the forward local rules to such a diagram since the
entries in the cells are not just $1$'s and $0$'s, and, even if they
should be (such as in the filling of the rectangle on the left of
Figure~\ref{fig:6}), there could be several $1$'s in a column or in a
row. To remedy this, we ``separate" the entries. 
That is, we construct a diagram with more rows and columns so
that entries which are originally 
in the same column or in the same row are put in different rows in
the larger diagram, and that an entry $m$ is replaced by $m$ $1$'s in
the new diagram, all of which placed in different rows and columns.
For this ``separation" we have two choices for the columns
and two choices for the rows: either we ``separate" entries in a row by
putting them into a chain from bottom/left to top/right, or we
``separate" them by putting them into a chain from top/left to
bottom/right, the same being true for entries in a column. 
In total, this gives $2\times2=4$ choices of separation.
As we shall see, two out of these four (the Second and Third Variant
below) are, in fact, equivalent
modulo a reflection of growth diagrams.
Interestingly, this
equivalence can only be seen directly from the growth diagram
algorithms, but not from the insertion algorithms that go with them.
(The First and Fourth Variant are also related, but in a much more
subtle way. See \cite[Sec.~3.2]{LeeuAH}.) 

\subsection{First variant: RSK}
The first variant (described in detail in \cite[Sec.~4.1]{RobyAA}) 
generalizes the Robinson--Schensted--Knuth (RSK)
correspondence. It is defined for arbitrary fillings of a
Ferrers shape with non-negative integers. 

Let us consider the filling of the $2\times 2$
square on the left of Figure~\ref{fig:6a}. (The labellings of the
corners of the cells should be ignored at this point.)
The filling is now converted into the $0$-$1$-filling of a larger
shape (where, again,
$1$'s are represented by X's and $0$'s are suppressed). 
If a cell is filled with entry $m$, 
we replace $m$ by a chain of $m$ X's arranged from bottom/left to
top/right. If there should be several entries in a column
then we arrange the chains coming from the entries of the column
as well from bottom/left to top/right. We do the same for the rows.
The diagram on the right of Figure~\ref{fig:6a} shows the result of
this conversion when applied to the filling on the left of
Figure~\ref{fig:6a}. (Still, the labellings of the corners of the
cells in the augmented diagram should be ignored at this point.)
In the figure, the original columns and rows are indicated by thick
lines, whereas the newly created columns and rows are indicated by
thin lines.

\begin{figure}[h]
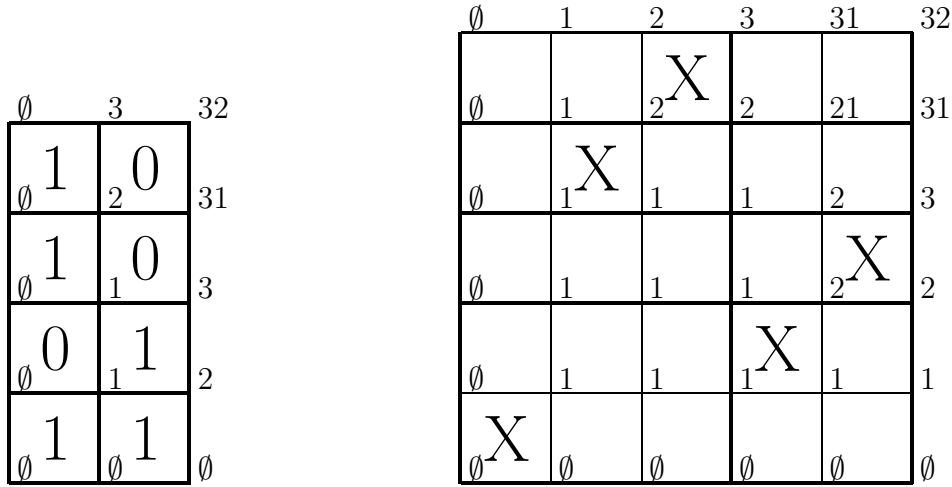

$$
\Einheit.4cm
\Pfad(0,12),111111\endPfad
\Pfad(0,9),111111\endPfad
\Pfad(0,6),111111\endPfad
\Pfad(0,3),111111\endPfad
\Pfad(0,0),111111\endPfad
\Pfad(0,0),222222222222\endPfad
\Pfad(3,0),222222222222\endPfad
\Pfad(6,0),222222222222\endPfad
\Label\ro{\text {\Huge $1$}}(1,1)
\Label\ro{\text {\Huge $1$}}(1,7)
\Label\ro{\text {\Huge $1$}}(1,10)
\Label\ro{\text {\Huge $1$}}(4,1)
\Label\ro{\text {\Huge $1$}}(4,4)
\Label\ro{\text {\Huge $0$}}(1,4)
\Label\ro{\text {\Huge $0$}}(4,7)
\Label\ro{\text {\Huge $0$}}(4,10)
\Label\ro{\emptyset}(0,0)
\Label\ro{\emptyset}(0,3)
\Label\ro{\emptyset}(0,6)
\Label\ro{\emptyset}(0,9)
\Label\ro{\emptyset}(0,12)
\Label\ro{\emptyset}(3,0)
\Label\ro{1}(3,3)
\Label\ro{1}(3,6)
\Label\ro{2}(3,9)
\Label\ro{3}(3,12)
\Label\ro{\emptyset}(6,0)
\Label\ro{2}(6,3)
\Label\ro{3}(6,6)
\Label\ro{\hphantom{1}31}(6,9)
\Label\ro{\hphantom{1}32}(6,12)
\hbox{\hskip6cm}
\Pfad(0,15),111111111111111\endPfad
\Pfad(0,12),111111111111111\endPfad
\Pfad(0,9),111111111111111\endPfad
\Pfad(0,6),111111111111111\endPfad
\PfadDicke{.3pt}
\Pfad(0,3),111111111111111\endPfad
\PfadDicke{1pt}
\Pfad(0,0),111111111111111\endPfad
\Pfad(15,0),222222222222222\endPfad
\PfadDicke{.3pt}
\Pfad(12,0),222222222222222\endPfad
\PfadDicke{1pt}
\Pfad(9,0),222222222222222\endPfad
\PfadDicke{.3pt}
\Pfad(6,0),222222222222222\endPfad
\Pfad(3,0),222222222222222\endPfad
\PfadDicke{1pt}
\Pfad(0,0),222222222222222\endPfad
\Label\ro{\emptyset}(0,0)
\Label\ro{\emptyset}(0,3)
\Label\ro{\emptyset}(0,6)
\Label\ro{\emptyset}(0,9)
\Label\ro{\emptyset}(0,12)
\Label\ro{\emptyset}(0,15)
\Label\ro{1}(3,15)
\Label\ro{1}(3,12)
\Label\ro{1}(3,9)
\Label\ro{1}(3,6)
\Label\ro{1}(3,3)
\Label\ro{\emptyset}(3,0)
\Label\ro{2}(6,15)
\Label\ro{2}(6,12)
\Label\ro{1}(6,9)
\Label\ro{1}(6,6)
\Label\ro{1}(6,3)
\Label\ro{\emptyset}(6,0)
\Label\ro{3}(9,15)
\Label\ro{2}(9,12)
\Label\ro{1}(9,9)
\Label\ro{1}(9,6)
\Label\ro{1}(9,3)
\Label\ro{\emptyset}(9,0)
\Label\ro{\hphantom{1}31}(12,15)
\Label\ro{\hphantom{1}21}(12,12)
\Label\ro{2}(12,9)
\Label\ro{2}(12,6)
\Label\ro{1}(12,3)
\Label\ro{\emptyset}(12,0)
\Label\ro{\hphantom{1}32}(15,15)
\Label\ro{\hphantom{1}31}(15,12)
\Label\ro{3}(15,9)
\Label\ro{2}(15,6)
\Label\ro{1}(15,3)
\Label\ro{\emptyset}(15,0)
\Label\ro{\text {\Huge X}}(1,1)
\Label\ro{\text {\Huge X}}(4,10)
\Label\ro{\text {\Huge X}}(7,13)
\Label\ro{\text {\Huge X}}(10,4)
\Label\ro{\text {\Huge X}}(13,7)
\hskip6cm
$$
\caption{Growth diagrams and RSK}
\label{fig:6}
\end{figure}

Now we can apply the forward rules (F1)--(F6) to the augmented
diagram. That is, we label all the corners of the cells on the left
side and the bottom side of the augmented diagram by $\emptyset$, and
then we apply (F1)--(F6) to determine the labels of all the other corners.
Subsequently, we ``shrink back" the augmented diagram,
that is, we record only the labels of the corners located at the
intersections of thick lines. This yields the labels on the left of
Figure~\ref{fig:6a}. For a different example see Figure~\ref{fig:6}.

The labellings by partitions that one obtains in this manner have again the 
property that a partition is contained in its right neighbour and in
its top neighbour. In addition, two neighbouring partitions
differ by a {\it horizontal strip}, that is, by a set of squares no
two of which are in the same column.

\begin{figure}[h]
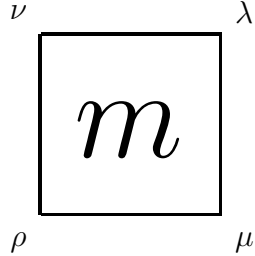

$$
\Pfad(0,0),1111222255556666\endPfad
\Label\lu{\rh}(0,0)
\Label\ru{\mu}(4,0)
\Label\lo{\nu}(0,4)
\Label\ro{\la}(4,4)
\thicklines
\Label\ro{\raise13pt\hbox{\hskip15pt\HUGE m}}(1,1)
\hskip2.4cm
$$
\vskip10pt
\caption{A cell filled with a non-negative integer $m$}
\label{fig:1a}
\end{figure}

Clearly, we could also define local forward and backward rules {\it
directly} on the original (smaller) diagram. 
This has been done in detail in \cite[Sec.~3.1]{LeeuAH}.
For the sake of completeness, we recall these rules here in a
slightly different, but of course equivalent, fashion. 
Consider the cell in
Figure~\ref{fig:1a}, filled by a non-negative integer $m$, and
labelled by the partitions $\rh$, $\mu$, $\nu$, where
$\rh\subseteq\mu$ and $\rh\subseteq\nu$, $\mu$ and $\rh$ differ by a
horizontal strip, and $\nu$ and $\rh$ differ by a horizontal strip.
Then $\la$ is determined by the following algorithm:
\begin{enumerate}
\item[(F$^1$0)]Set $\text {CARRY}:=m$ and $i:=1$.
\item[(F$^1$1)]Set $\la_i:=\max\{\mu_i,\nu_i\}+\text {CARRY}$.
\item[(F$^1$2)]If $\la_i=0$, then stop. The output of the algorithm
is $\la=(\la_1,\la_2,\dots,\la_{i-1})$. If not, then
set $\text {CARRY}:=\min\{\mu_i,\nu_i\}-\rh_i$ and $i:=i+1$.
Go to (F$^1$1).
\end{enumerate}

Conversely,
given $\mu,\nu,\la$, where $\mu\subseteq\la$ and $\nu\subseteq\la$,
where $\la$ and $\mu$ differ by a horizontal strip, and 
where $\la$ and $\nu$ differ by a horizontal strip, 
the backward algorithm works in the following way:
\begin{enumerate}
\item[(B$^1$0)]Set $i:=\max\{j:\la_j\text{ is positive}\}$, and $\text
{CARRY}:=0$.
\item[(B$^1$1)]Set $\rh_i:=\min\{\mu_i,\nu_i\}-\text{CARRY}$.
\item[(B$^1$2)]Set $\text {CARRY}:=\la_i-\max\{\mu_i,\nu_i\}$ and $i:=i-1$.
If $i=0$, then stop. The output of the algorithm is
$\rh=(\rh_1,\rh_2,\dots)$ and $m=\text {CARRY}$. If not, go to (B$^1$1).
\end{enumerate}

The extension of Theorem~\ref{thm:1} then reads as follows.
In the statement of the theorem below (and also later), 
for two partitions with $\nu\subseteq\mu$ we write $\mu/\nu$
for the difference of the Ferrers diagrams corresponding to $\mu$ and
$\nu$, respectively, that is, for the set of squares which belong to
the Ferrers diagram of $\mu$ but not to the Ferrers diagram of $\nu$.

\begin{Theorem} \label{thm:1a}
Let $F$ be a Ferrers shape given by the $D$-$R$-sequence
$w=w_1w_2\dots w_k$.
Fillings of $F$ with non-negative integers 
are in bijection with sequences
$(\emptyset=\la^0,\la^1,\dots,\la^{k}=\emptyset)$, where
$\la^i/\la^{i-1}$ is a horizontal strip if $w_i=R$, whereas
$\la^{i-1}/\la^i$ is a horizontal strip if $w_i=D$.
\end{Theorem}

It is a now well-known fact (cf.\ \cite[Sec.~4.1]{RobyAA}) that,
in the case that the Ferrers shape is a rectangle,
the bijection in Theorem~\ref{thm:1a} is equivalent to the
{\it Robinson--Schensted--Knuth correspondence}.
Namely, assuming that $F$ is a $p\times q$ rectangle,
according to Theorem~\ref{thm:1a},
the sequences $(\la^0,\la^1,\dots,\la^{p+q})$ 
which we read off the
top side and the right side of the square are sequences with
$\emptyset=\la^0\subseteq\la^1\subseteq\dots\subseteq\la^q
\supseteq\dots\supseteq\la^{p+q-1}\supseteq\la^{p+q}=\emptyset$, 
where $\la^i/\la^{i-1}$ is a horizontal strip for
$i=1,2,\dots,q$, and where $\la^{i}/\la^{i+1}$ is a horizontal strip for
$i=q,q+1,\dots,p+q-1$.
In their turn, these are in bijection with pairs $(P,Q)$ of semistandard
tableaux of the same shape, the common shape consisting of
$n$ squares, where $n$ is the sum of all the entries of the
filling. (See \cite{RobyAA}. The semistandard tableau $Q$ is defined
by the first half of the sequence, 
$\emptyset=\la^0\subseteq\la^1\subseteq\dots\subseteq\la^q$,
the entry $i$ being put in the squares by which $\la^i$ and $\la^{i-1}$
differ, and, similarly, the semistandard tableau $P$ is defined
by the second half of the sequence, 
$\emptyset=\la^{p+q}\subseteq\la^{p+q-1}\subseteq\dots\subseteq\la^q$.)
It is then a theorem 
that the bijection between fillings and pairs of semistandard tableaux
defined by the growth diagram on the rectangle coincides with
the Robinson--Schensted--Knuth correspondence, $P$ being the {\it insertion
tableau}, and $Q$ being the {\it recording tableau}.
(The reader is referred to \cite[Ch.~4]{FultAC},
\cite{KnutAA},
\cite[Sec.~4.8]{SagaAQ}, \cite[Sec.~7.11]{StanBI} for extensive
information on the Robinson--Schensted--Knuth correspondence.)

Figure~\ref{fig:6} is meant to illustrate this, using an
example that will serve as a running example here as well as 
for the other three variants. Let us recall 
that the RSK correspondence starts with a rectangular filling
as on the left of Figure~\ref{fig:6}, transforms the filling
into a two-rowed array, and then transforms the two-rowed array
into a pair of semistandard tableaux by an insertion procedure.
The two-rowed array is obtained from the filling, by considering
the entry, $m$ say, in the $i$-th row (from below) and $j$-th
column (from the left), and recording $m$ pairs $\binom ji$.
These pairs are then ordered into a two-rowed array such that
the entries in the top row are weakly increasing, and,
in the bottom row, entries must be weakly increasing 
below equal entries in the top row. Thus, the filling
in Figure~\ref{fig:6a} corresponds to the two-rowed array
$$\begin{pmatrix} 1\ 1\ 1\ 2\ 2\\1\ 2\ 2\ 1\ 1\end{pmatrix},$$
while the filling
in Figure~\ref{fig:6} corresponds to the two-rowed array
\begin{equation} \label{eq:2r} 
\begin{pmatrix} 1\ 1\ 1\ 2\ 2\\1\ 3\ 4\ 1\ 2\end{pmatrix}.
\end{equation}
Now the bottom entries are inserted according to {\it row insertion},
in which an element bumps the first entry in a row which is strictly
larger.
The top entries keep track of where the 
individual insertions stop. If this is applied to the two-rowed 
array in \eqref{eq:2r}, one obtains
$$(P,Q)=\left(\begin{matrix} 1\ 1\ 2\\3\ 4\ \hphantom{2}\end{matrix}\ ,\
\begin{matrix} 1\ 1\ 1\\2\ 2\ \hphantom{2}\end{matrix}\right),$$
which is indeed in agreement with the increasing sequences
of partitions along the right side and the top side of the
rectangle, respectively.

\medskip
Again, in addition to its local description,
the bijection in Theorem~\ref{thm:1a} has also a {\it global\/}
description. It is again a consequence of Greene's theorem (stated here as
Theorem~\ref{thm:3}) and the description of the bijection based on
``separation" of entries along columns and rows.
In order to formulate the
result, we adapt our previous definitions: a {\it NE-chain} of a
filling is a sequence of non-zero entries in the filling such that any entry
in the sequence is weakly above and weakly to the right of the
preceding entry in the sequence. The {\it length\/} of such a
NE-chain is defined as the {\it sum} of all the entries in the chain.
On the other hand, a {\it se-chain} of a
filling is a sequence of non-zero entries in the filling such that any entry
in the sequence is strictly below and strictly to the right of the
preceding entry in the sequence. In contrast
to NE-chains, we define
the {\it length\/} of a se-chain as the {\it number} of entries in the
chain. (These definitions can be best motivated as weak, respectively
strict, chains of balls in the matrix--ball model for
fillings described in \cite[Sec.~4.2]{FultAC}.)

\begin{Theorem} \label{thm:3a}
Given a diagram with empty partitions labelling all the
corners along the left side and the bottom side of the Ferrers shape,
which has been completed according to RSK,
the partition $\la=(\la_1,\la_2,\dots,\la_\ell)$ labelling corner $c$ 
satisfies the following two properties:
\begin{enumerate}
\item [(G$^1$1)]For any $k$, the maximal sum of all the entries
in a collection of $k$ pairwise disjoint
NE-chains situated in the rectangular region to the left and
below of $c$ is equal to $\la_1+\la_2+\dots+\la_k$.
\item [(G$^1$2)]Fix a positive integer $k$ and consider collections
of $k$ se-chains with the property that no entry $e$ can be in more
than $e$ of these se-chains. Then
the maximal cardinality of the multiset\/ {\em(!)} union of the se-chains 
in such a collection of se-chains
situated in the rectangular region to the left and
below of $c$ is equal to $\la'_1+\la'_2+\dots+\la'_k$, where $\la'$
denotes the partition conjugate to $\la$.
\end{enumerate}
In particular, $\la_1$ is the length of the longest NE-chain
in the rectangular region to the left and below of $c$, and $\la'_1$
is the length of the longest se-chain in the same
rectangular region.
\end{Theorem}

\subsection{Second variant: dual RSK}
The second variant
generalizes the ``dual correspondence" of Knuth \cite{KnutAA},
which we abbreviate as {\it dual RSK}. 
It is only defined for $0$-$1$-fillings of a
Ferrers shape. (In contrast to Sections~\ref{sec:growth} and
\ref{sec:part}, these can, however, be arbitrary, that is,
there can be several $1$'s in a row or a column.)

\begin{figure}[h]
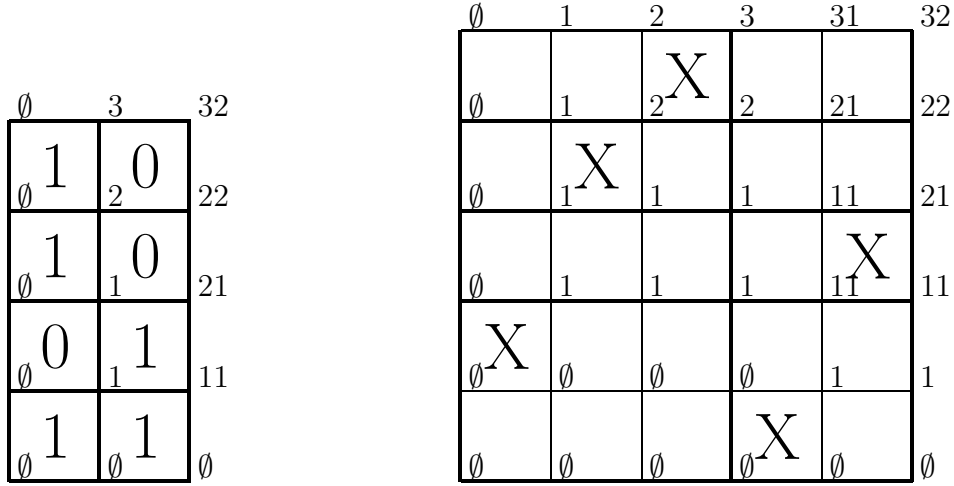

$$
\Einheit.4cm
\Pfad(0,12),111111\endPfad
\Pfad(0,9),111111\endPfad
\Pfad(0,6),111111\endPfad
\Pfad(0,3),111111\endPfad
\Pfad(0,0),111111\endPfad
\Pfad(0,0),222222222222\endPfad
\Pfad(3,0),222222222222\endPfad
\Pfad(6,0),222222222222\endPfad
\Label\ro{\text {\Huge $1$}}(1,1)
\Label\ro{\text {\Huge $1$}}(1,7)
\Label\ro{\text {\Huge $1$}}(1,10)
\Label\ro{\text {\Huge $1$}}(4,1)
\Label\ro{\text {\Huge $1$}}(4,4)
\Label\ro{\text {\Huge $0$}}(1,4)
\Label\ro{\text {\Huge $0$}}(4,7)
\Label\ro{\text {\Huge $0$}}(4,10)
\Label\ro{\emptyset}(0,0)
\Label\ro{\emptyset}(0,3)
\Label\ro{\emptyset}(0,6)
\Label\ro{\emptyset}(0,9)
\Label\ro{\emptyset}(0,12)
\Label\ro{\emptyset}(3,0)
\Label\ro{1}(3,3)
\Label\ro{1}(3,6)
\Label\ro{2}(3,9)
\Label\ro{3}(3,12)
\Label\ro{\emptyset}(6,0)
\Label\ro{\hphantom{1}11}(6,3)
\Label\ro{\hphantom{1}21}(6,6)
\Label\ro{\hphantom{1}22}(6,9)
\Label\ro{\hphantom{1}32}(6,12)
\hbox{\hskip6cm}
\Pfad(0,15),111111111111111\endPfad
\Pfad(0,12),111111111111111\endPfad
\Pfad(0,9),111111111111111\endPfad
\Pfad(0,6),111111111111111\endPfad
\PfadDicke{.3pt}
\Pfad(0,3),111111111111111\endPfad
\PfadDicke{1pt}
\Pfad(0,0),111111111111111\endPfad
\Pfad(15,0),222222222222222\endPfad
\PfadDicke{.3pt}
\Pfad(12,0),222222222222222\endPfad
\PfadDicke{1pt}
\Pfad(9,0),222222222222222\endPfad
\PfadDicke{.3pt}
\Pfad(6,0),222222222222222\endPfad
\Pfad(3,0),222222222222222\endPfad
\PfadDicke{1pt}
\Pfad(0,0),222222222222222\endPfad
\Label\ro{\emptyset}(0,0)
\Label\ro{\emptyset}(0,3)
\Label\ro{\emptyset}(0,6)
\Label\ro{\emptyset}(0,9)
\Label\ro{\emptyset}(0,12)
\Label\ro{\emptyset}(0,15)
\Label\ro{1}(3,15)
\Label\ro{1}(3,12)
\Label\ro{1}(3,9)
\Label\ro{1}(3,6)
\Label\ro{\emptyset}(3,3)
\Label\ro{\emptyset}(3,0)
\Label\ro{2}(6,15)
\Label\ro{2}(6,12)
\Label\ro{1}(6,9)
\Label\ro{1}(6,6)
\Label\ro{\emptyset}(6,3)
\Label\ro{\emptyset}(6,0)
\Label\ro{3}(9,15)
\Label\ro{2}(9,12)
\Label\ro{1}(9,9)
\Label\ro{1}(9,6)
\Label\ro{\emptyset}(9,3)
\Label\ro{\emptyset}(9,0)
\Label\ro{\hphantom{1}31}(12,15)
\Label\ro{\hphantom{1}21}(12,12)
\Label\ro{\hphantom{1}11}(12,9)
\Label\ro{\hphantom{1}11}(12,6)
\Label\ro{1}(12,3)
\Label\ro{\emptyset}(12,0)
\Label\ro{\hphantom{1}32}(15,15)
\Label\ro{\hphantom{1}22}(15,12)
\Label\ro{\hphantom{1}21}(15,9)
\Label\ro{\hphantom{1}11}(15,6)
\Label\ro{1}(15,3)
\Label\ro{\emptyset}(15,0)
\Label\ro{\text {\Huge X}}(1,4)
\Label\ro{\text {\Huge X}}(4,10)
\Label\ro{\text {\Huge X}}(7,13)
\Label\ro{\text {\Huge X}}(10,1)
\Label\ro{\text {\Huge X}}(13,7)
\hskip6cm
$$
\caption{Growth diagrams and dual RSK}
\label{fig:8}
\end{figure}

As an illustration for this second variant, we use again
the filling on the left of Figure~\ref{fig:6}, which is
reproduced on the left of Figure~\ref{fig:8}. 
(The labellings of the
corners of the cells should be ignored at this point.)
The filling is now converted into a $0$-$1$-filling of a larger
shape, but in a different way as before (where, again,
$1$'s will be represented by X's and $0$'s are suppressed). 
Namely, if there should be several $1$'s in a column
then, as in the first variant, 
we arrange them from bottom/left to top/right. However,
we do the ``opposite" for the rows, that is,
if there should be several $1$'s in a row
then we arrange them from top/left to bottom/right.
The diagram on the right of Figure~\ref{fig:8} shows the result of
this conversion when applied to the filling on the left of
Figure~\ref{fig:8}. 

Now, as in the first variant, 
we apply the forward rules (F1)--(F6) to the augmented
diagram, and, once this is done,
we ``shrink back" the augmented diagram, as before.
This yields the labels on the left of
Figure~\ref{fig:8}. 

Of course,
the labellings by partitions that one obtains in this manner have again the 
property that a partition is contained in its right neighbour and in
its top neighbour. However, now two partitions which are
horizontal neighbours 
differ by a {\it horizontal strip}, whereas two partitions
which are {\it vertical\/} neighbours differ by a
{\it vertical\/} strip, 
that is, by a set of squares no
two of which are in the same row.

{\it Direct} local forward and backward rules
on the original (smaller) diagram are also available for this
variant, see \cite[Sec.~6]{LeeuAH}. Our presentation is again
slightly different, but equivalent.
Namely, consider the cell in
Figure~\ref{fig:1a}, filled by $m=0$ or $m=1$ and
labelled by the partitions $\rh$, $\mu$, $\nu$, where
$\rh\subseteq\mu$ and $\rh\subseteq\nu$, $\mu$ and $\rh$ differ by a
horizontal strip, and $\nu$ and $\rh$ differ by a vertical strip.
Then $\la$ is determined by the following algorithm:
\begin{enumerate}
\item[(F$^2$0)]Set $\text {CARRY}:=m$ and $i:=1$.
\item[(F$^2$1)]Set $\la_i:=\max\{\mu_i+\text {CARRY},\nu_i\}$.
\item[(F$^2$2)]If $\la_i=0$, then stop. The output of the algorithm
is $\la=(\la_1,\la_2,\dots,\la_{i-1})$. If not, then
set $\text {CARRY}:=\min\{\mu_i+\text {CARRY},\nu_i\}-\rh_i$ 
and $i:=i+1$.
Go to (F$^2$1).
\end{enumerate}

Conversely,
given $\mu,\nu,\la$, where $\mu\subseteq\la$ and $\nu\subseteq\la$,
where $\la$ and $\mu$ differ by a vertical strip, and 
where $\la$ and $\nu$ differ by a horizontal strip, 
the backward algorithm works in the following way:
\begin{enumerate}
\item[(B$^2$0)]Set $i:=\max\{j:\la_j\text{ is positive}\}$ and $\text
{CARRY}:=0$.
\item[(B$^2$1)]Set $\rh_i:=\min\{\mu_i,\nu_i-\text {CARRY}\}$.
\item[(B$^2$2)]Set $\text {CARRY}:=\la_i-\max\{\mu_i,\nu_i-\text {CARRY}\}$ 
and $i:=i-1$.
If $i=0$, then stop. The output of the algorithm is
$\rh=(\rh_1,\rh_2,\dots)$ and $m=\text {CARRY}$. If not, go to (B$^2$1).
\end{enumerate}

The corresponding extension of Theorem~\ref{thm:1} then reads as follows.

\begin{Theorem} \label{thm:1b}
Let $F$ be a Ferrers shape given by the $D$-$R$-sequence
$w=w_1w_2\dots w_k$.
Then $0$-$1$-fillings of $F$ with non-negative integers 
are in bijection with sequences
$(\emptyset=\la^0,\la^1,\dots,\la^{k}=\emptyset)$, where
$\la^i/\la^{i-1}$ is a horizontal strip if $w_i=R$, whereas
$\la^{i-1}/\la^i$ is a vertical strip if $w_i=D$.
\end{Theorem}

In the case that the Ferrers shape is a rectangle,
the bijection in Theorem~\ref{thm:1b} is equivalent to 
dual RSK.
Namely, assuming that $F$ is a $p\times q$ rectangle,
according to Theorem~\ref{thm:1b},
the sequences $(\la^0,\la^1,\dots,\la^{p+q})$ 
which we read off the
top side and the right side of the square are sequences with
$\emptyset=\la^0\subseteq\la^1\subseteq\dots\subseteq\la^q
\supseteq\dots\supseteq\la^{p+q-1}\supseteq\la^{p+q}=\emptyset$, 
where $\la^i/\la^{i-1}$ is a horizontal strip for
$i=1,2,\dots,q$, and where $\la^{i}/\la^{i+1}$ is a vertical strip for
$i=q,q+1,\dots,p+q-1$.
In their turn, these are in bijection with pairs $(P,Q)$,
where $Q$ and the transpose of $P$ are semistandard
tableaux of the same shape, the common shape consisting of
$n$ squares, where $n$ is the sum of all the entries of the
filling. The semistandard tableau $Q$ is defined
by the first half of the sequence, 
$\emptyset=\la^0\subseteq\la^1\subseteq\dots\subseteq\la^q$,
as before, while $P$ is defined
by the second half of the sequence, 
$\emptyset=\la^{p+q}\subseteq\la^{p+q-1}\subseteq\dots\subseteq\la^q$,
as before. Since, in the latter chain of partitions, successive
partitions differ by {\it vertical\/} strips, $P$ is not a
semistandard tableau, but its transpose is.
This bijection between fillings and pairs $(P,Q)$ coincides with
dual RSK.
(We refer the reader to \cite[App.~A.4.3]{FultAC},
\cite[Sec.~5]{KnutAA},
\cite[Sec.~4.8]{SagaAQ}, \cite[Sec.~7.14]{StanBI} for more
information on the dual correspondence.)

Figure~\ref{fig:8} illustrates this with our running example.
Let us recall 
that dual RSK starts with a rectangular filling, consisting of $0$'s
and $1$'s,
as on the left of Figure~\ref{fig:8}.
The filling is now transformed
into a two-rowed array, by recording a pair $\binom ji$
for a $1$ in the $i$-th row (from below) and $j$-th
column (from the left). Subseqently, 
the pairs are ordered into a two-rowed array as before, so that we
obtain again the two-rowed array \eqref{eq:2r}.
Now the bottom entries are inserted according to {\it column
insertion} (cf.\ \cite[Sec.~5, Algorithm INSERT*; the tableaux 
there must be transposed to obtain our version here]{KnutAA}), 
in which an element bumps the
first entry in a column which is larger than or equal to it.
Again, the top entries keep track of where the 
individual insertions stop. If this is applied to the two-rowed 
array in \eqref{eq:2r}, one obtains
\begin{equation} \label{eq:P,Q}
(P^t,Q^t)=\left(\begin{matrix} 1\ 1\\2\ 3\\4\ \hphantom{2}\end{matrix}\ ,\
\begin{matrix} 1\ 2\\1\ 2\\1\ \hphantom{2}\end{matrix}\right),
\end{equation}
which is indeed in agreement with the increasing sequences
of partitions along the right side and the top side of the
rectangle, respectively, if one transposes both arrays in
\eqref{eq:P,Q}.

\medskip
Again, in addition to its local description,
the bijection in Theorem~\ref{thm:1b} has also a {\it global\/}
description. It is again a consequence of Greene's theorem (stated here as
Theorem~\ref{thm:3}) and the description of the bijection based on
``separation" of entries along columns and rows.
In order to formulate the global description of the bijection in
terms of increasing and decreasing chains, we need to define
a {\it nE-chain} of a
filling to be a sequence of $1$'s in the filling such that any $1$
in the sequence is strictly above and weakly to the right of the
preceding $1$ in the sequence. 
Furthermore, we define a {\it Se-chain} of a
filling to be a sequence of $1$'s in the filling such that any $1$ in the
sequence is weakly below and strictly to the right of the 
preceding $1$ in the sequence. 
The {\it length\/} of a nE-chain or a Se-chain
is defined as the number of $1$'s in the chain.
(Again, these definitions can be best motivated as 
chains of balls in the matrix--ball model for
fillings described in \cite[Sec.~4.2]{FultAC}.)

\begin{Theorem} \label{thm:3b}
Given a diagram with empty partitions labelling all the
corners along the left side and the bottom side of the Ferrers shape,
which has been completed according to dual RSK,
the partition $\la=(\la_1,\la_2,\dots,\la_\ell)$ labelling corner $c$ 
satisfies the following two properties:
\begin{enumerate}
\item [(G$^2$1)]For any $k$, the maximal cardinality of the union of $k$
nE-chains situated in the rectangular region to the left and
below of $c$ is equal to $\la_1+\la_2+\dots+\la_k$.
\item [(G$^2$2)]For any $k$, the maximal cardinality of the union of $k$
Se-chains situated in the rectangular region to the left and
below of $c$ is equal to $\la'_1+\la'_2+\dots+\la'_k$, where $\la'$
denotes the partition conjugate to $\la$.
\end{enumerate}
In particular, $\la_1$ is the length of the longest nE-chain
in the rectangular region to the left and below of $c$, and $\la'_1$
is the length of the longest Se-chain in the same
rectangular region.
\end{Theorem}

\subsection{Third variant: RSK'}\label{sec:RSK'}
\begin{figure}[h]
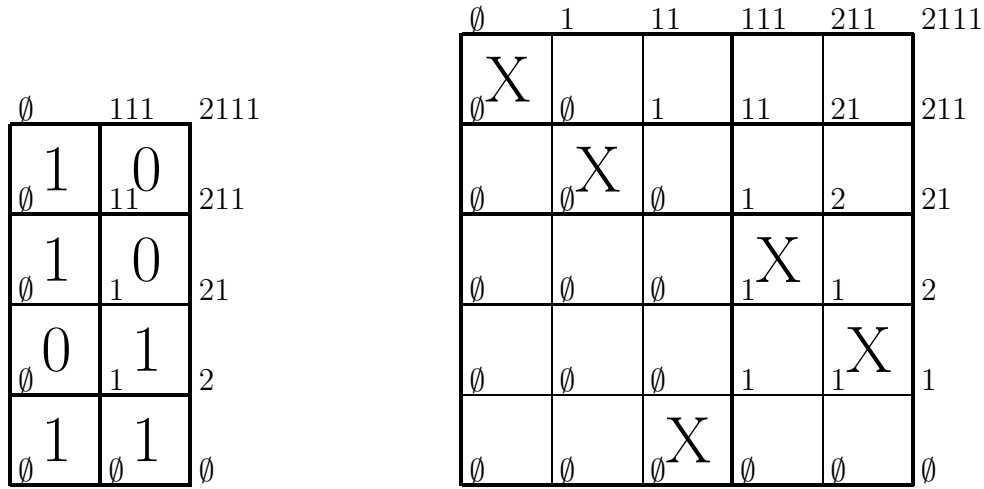

$$
\Einheit.4cm
\Pfad(0,12),111111\endPfad
\Pfad(0,9),111111\endPfad
\Pfad(0,6),111111\endPfad
\Pfad(0,3),111111\endPfad
\Pfad(0,0),111111\endPfad
\Pfad(0,0),222222222222\endPfad
\Pfad(3,0),222222222222\endPfad
\Pfad(6,0),222222222222\endPfad
\Label\ro{\text {\Huge $1$}}(1,1)
\Label\ro{\text {\Huge $1$}}(1,7)
\Label\ro{\text {\Huge $1$}}(1,10)
\Label\ro{\text {\Huge $1$}}(4,1)
\Label\ro{\text {\Huge $1$}}(4,4)
\Label\ro{\text {\Huge $0$}}(1,4)
\Label\ro{\text {\Huge $0$}}(4,7)
\Label\ro{\text {\Huge $0$}}(4,10)
\Label\ro{\emptyset}(0,0)
\Label\ro{\emptyset}(0,3)
\Label\ro{\emptyset}(0,6)
\Label\ro{\emptyset}(0,9)
\Label\ro{\emptyset}(0,12)
\Label\ro{\emptyset}(3,0)
\Label\ro{1}(3,3)
\Label\ro{1}(3,6)
\Label\ro{\hphantom{1}11}(3,9)
\Label\ro{\hphantom{11}111}(3,12)
\Label\ro{\emptyset}(6,0)
\Label\ro{2}(6,3)
\Label\ro{\hphantom{1}21}(6,6)
\Label\ro{\hphantom{11}211}(6,9)
\Label\ro{\hphantom{111}2111}(6,12)
\hbox{\hskip6cm}
\Pfad(0,15),111111111111111\endPfad
\Pfad(0,12),111111111111111\endPfad
\Pfad(0,9),111111111111111\endPfad
\Pfad(0,6),111111111111111\endPfad
\PfadDicke{.3pt}
\Pfad(0,3),111111111111111\endPfad
\PfadDicke{1pt}
\Pfad(0,0),111111111111111\endPfad
\Pfad(15,0),222222222222222\endPfad
\PfadDicke{.3pt}
\Pfad(12,0),222222222222222\endPfad
\PfadDicke{1pt}
\Pfad(9,0),222222222222222\endPfad
\PfadDicke{.3pt}
\Pfad(6,0),222222222222222\endPfad
\Pfad(3,0),222222222222222\endPfad
\PfadDicke{1pt}
\Pfad(0,0),222222222222222\endPfad
\Label\ro{\emptyset}(0,0)
\Label\ro{\emptyset}(0,3)
\Label\ro{\emptyset}(0,6)
\Label\ro{\emptyset}(0,9)
\Label\ro{\emptyset}(0,12)
\Label\ro{\emptyset}(0,15)
\Label\ro{1}(3,15)
\Label\ro{\emptyset}(3,12)
\Label\ro{\emptyset}(3,9)
\Label\ro{\emptyset}(3,6)
\Label\ro{\emptyset}(3,3)
\Label\ro{\emptyset}(3,0)
\Label\ro{\hphantom{1}11}(6,15)
\Label\ro{1}(6,12)
\Label\ro{\emptyset}(6,9)
\Label\ro{\emptyset}(6,6)
\Label\ro{\emptyset}(6,3)
\Label\ro{\emptyset}(6,0)
\Label\ro{\hphantom{11}111}(9,15)
\Label\ro{\hphantom{1}11}(9,12)
\Label\ro{1}(9,9)
\Label\ro{1}(9,6)
\Label\ro{1}(9,3)
\Label\ro{\emptyset}(9,0)
\Label\ro{\hphantom{11}211}(12,15)
\Label\ro{\hphantom{1}21}(12,12)
\Label\ro{2}(12,9)
\Label\ro{1}(12,6)
\Label\ro{1}(12,3)
\Label\ro{\emptyset}(12,0)
\Label\ro{\hphantom{111}2111}(15,15)
\Label\ro{\hphantom{11}211}(15,12)
\Label\ro{\hphantom{1}21}(15,9)
\Label\ro{2}(15,6)
\Label\ro{1}(15,3)
\Label\ro{\emptyset}(15,0)
\Label\ro{\text {\Huge X}}(1,13)
\Label\ro{\text {\Huge X}}(4,10)
\Label\ro{\text {\Huge X}}(7,1)
\Label\ro{\text {\Huge X}}(10,7)
\Label\ro{\text {\Huge X}}(13,4)
\hskip6cm
$$
\caption{Growth diagrams and RSK'}
\label{fig:7}
\end{figure}

The third variant which we have in mind is, in growth diagram terms,
the ``reflection" of the second variant. By ``reflection," we mean
reflection of growth diagrams in a diagonal. More precisely, 
given a $0$-$1$-filling of a Ferrers diagram $F$, we separate $1$'s
in the same column by arranging
them from top/left to bottom/right, while we separate $1$'s
in the same row by arranging them from bottom/left to top/right.
What we obtain when we apply this to our running example,
is shown in Figure~\ref{fig:7}.

We abbreviate this algorithm as {\it RSK'}.
Since, as we said, in growth diagram terms, this is just the second
variant, but reflected in a diagonal, 
the relevant facts have already been told when discussing dual RSK,
with one exception: we have to explain to which
insertion algorithm RSK' is equivalent.

In order to do so, we transform again our $0$-$1$-filling in
Figure~\ref{fig:7} into a two-rowed array. This is again done by
constructing pairs from the $1$'s in the filling, as before. 
However, the pairs
are now ordered in a different way. Namely, we order the pairs such
that the entries in the top row are weakly increasing, 
in the bottom row, however, entries must be {\it decreasing}
below equal entries in the top row. Thus, the filling
in Figure~\ref{fig:7} corresponds to the two-rowed array
\begin{equation} \label{eq:2rc} 
\begin{pmatrix} 1\ 1\ 1\ 2\ 2\\4\ 3\ 1\ 2\ 1\end{pmatrix}.
\end{equation}
Now we apply row insertion to the bottom entries to construct a
semistandard tableau $P$ and use the top entries to record
the insertions in the array $Q$. In general, $Q$ will
not be a semistandard tableau, but the transpose of $Q$ will be.
If we apply this procedure to our two-rowed array in \eqref{eq:2rc},
then we obtain the pair
$$(P,Q)=\left(\begin{matrix} 1\ 1\\2\ \hphantom{2}\\3\ \hphantom{2}
\\4\ \hphantom{2}\end{matrix}\ ,\
\begin{matrix} 1\ 2\\1\ \hphantom{2}\\1\ \hphantom{2}\\2\ 
\hphantom{2}\end{matrix}\right),$$
which is indeed in agreement with the increasing sequences
of partitions along the right side and the top side of the
rectangle in Figure~\ref{fig:7}, respectively.
(An equivalent insertion algorithm is described in
\cite[App.~A.4.3, (1d), (2d)]{FultAC}.)

\subsection{Fourth variant: dual RSK'}
Our last variant is, as the first variant, 
defined for arbitrary fillings of a
Ferrers shape with non-negative integers. We abbreviate it by
{\it dual RSK'}.

\begin{figure}[h]
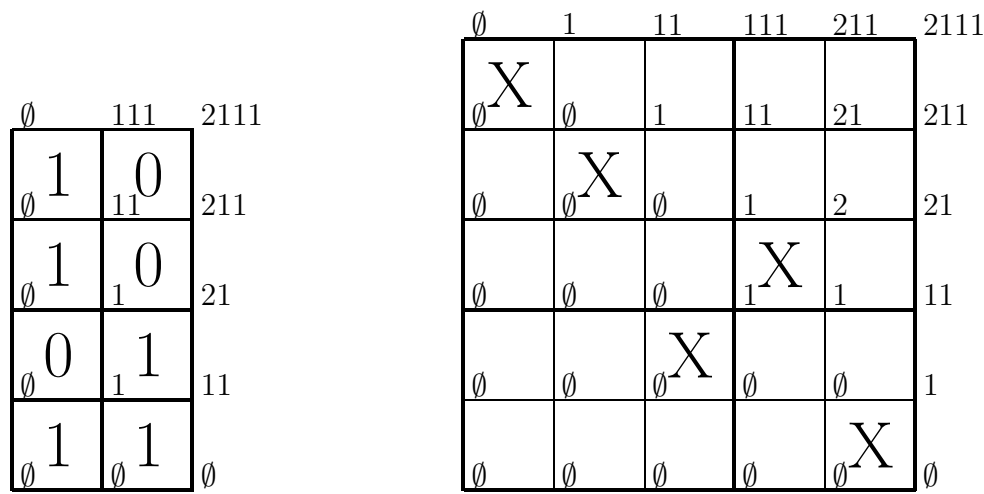

$$
\Einheit.4cm
\Pfad(0,12),111111\endPfad
\Pfad(0,9),111111\endPfad
\Pfad(0,6),111111\endPfad
\Pfad(0,3),111111\endPfad
\Pfad(0,0),111111\endPfad
\Pfad(0,0),222222222222\endPfad
\Pfad(3,0),222222222222\endPfad
\Pfad(6,0),222222222222\endPfad
\Label\ro{\text {\Huge $1$}}(1,1)
\Label\ro{\text {\Huge $1$}}(1,7)
\Label\ro{\text {\Huge $1$}}(1,10)
\Label\ro{\text {\Huge $1$}}(4,1)
\Label\ro{\text {\Huge $1$}}(4,4)
\Label\ro{\text {\Huge $0$}}(1,4)
\Label\ro{\text {\Huge $0$}}(4,7)
\Label\ro{\text {\Huge $0$}}(4,10)
\Label\ro{\emptyset}(0,0)
\Label\ro{\emptyset}(0,3)
\Label\ro{\emptyset}(0,6)
\Label\ro{\emptyset}(0,9)
\Label\ro{\emptyset}(0,12)
\Label\ro{\emptyset}(3,0)
\Label\ro{1}(3,3)
\Label\ro{1}(3,6)
\Label\ro{\hphantom{1}11}(3,9)
\Label\ro{\hphantom{11}111}(3,12)
\Label\ro{\emptyset}(6,0)
\Label\ro{\hphantom{1}11}(6,3)
\Label\ro{\hphantom{1}21}(6,6)
\Label\ro{\hphantom{11}211}(6,9)
\Label\ro{\hphantom{111}2111}(6,12)
\hbox{\hskip6cm}
\Pfad(0,15),111111111111111\endPfad
\Pfad(0,12),111111111111111\endPfad
\Pfad(0,9),111111111111111\endPfad
\Pfad(0,6),111111111111111\endPfad
\PfadDicke{.3pt}
\Pfad(0,3),111111111111111\endPfad
\PfadDicke{1pt}
\Pfad(0,0),111111111111111\endPfad
\Pfad(15,0),222222222222222\endPfad
\PfadDicke{.3pt}
\Pfad(12,0),222222222222222\endPfad
\PfadDicke{1pt}
\Pfad(9,0),222222222222222\endPfad
\PfadDicke{.3pt}
\Pfad(6,0),222222222222222\endPfad
\Pfad(3,0),222222222222222\endPfad
\PfadDicke{1pt}
\Pfad(0,0),222222222222222\endPfad
\Label\ro{\emptyset}(0,0)
\Label\ro{\emptyset}(0,3)
\Label\ro{\emptyset}(0,6)
\Label\ro{\emptyset}(0,9)
\Label\ro{\emptyset}(0,12)
\Label\ro{\emptyset}(0,15)
\Label\ro{1}(3,15)
\Label\ro{\emptyset}(3,12)
\Label\ro{\emptyset}(3,9)
\Label\ro{\emptyset}(3,6)
\Label\ro{\emptyset}(3,3)
\Label\ro{\emptyset}(3,0)
\Label\ro{\hphantom{1}11}(6,15)
\Label\ro{1}(6,12)
\Label\ro{\emptyset}(6,9)
\Label\ro{\emptyset}(6,6)
\Label\ro{\emptyset}(6,3)
\Label\ro{\emptyset}(6,0)
\Label\ro{\hphantom{11}111}(9,15)
\Label\ro{\hphantom{1}11}(9,12)
\Label\ro{1}(9,9)
\Label\ro{1}(9,6)
\Label\ro{\emptyset}(9,3)
\Label\ro{\emptyset}(9,0)
\Label\ro{\hphantom{11}211}(12,15)
\Label\ro{\hphantom{1}21}(12,12)
\Label\ro{2}(12,9)
\Label\ro{1}(12,6)
\Label\ro{\emptyset}(12,3)
\Label\ro{\emptyset}(12,0)
\Label\ro{\hphantom{111}2111}(15,15)
\Label\ro{\hphantom{11}211}(15,12)
\Label\ro{\hphantom{1}21}(15,9)
\Label\ro{\hphantom{1}11}(15,6)
\Label\ro{1}(15,3)
\Label\ro{\emptyset}(15,0)
\Label\ro{\text {\Huge X}}(1,13)
\Label\ro{\text {\Huge X}}(4,10)
\Label\ro{\text {\Huge X}}(7,4)
\Label\ro{\text {\Huge X}}(10,7)
\Label\ro{\text {\Huge X}}(13,1)
\hskip6cm
$$
\caption{Growth diagrams and dual RSK'}
\label{fig:9}
\end{figure}

Our running example serves once more to illustrate this fourth variant, 
see Figure~\ref{fig:9}.
Here, we separate entries in the following way.
If a cell is filled with entry $m$, 
we replace $m$ by a chain of $m$ X's arranged from top/left to
bottom/right. If there should be several entries in a column
then we arrange the chains coming from the entries of the column
as well from top/left to bottom/right. We do the same for the rows.
(In brief, everything is reversed compared to the first variant.)

Finally, as always,
we apply the forward rules (F1)--(F6) to the augmented
diagram, and, once this is done,
we ``shrink back" the augmented diagram. See Figure~\ref{fig:9}.

The labellings by partitions that one obtains in this manner have again the 
property that a partition is contained in its right neighbour and in
its top neighbour. However, here, two neighbouring partitions
differ by a {\it vertical strip}.

The corresponding {\it direct} local forward and backward rules
have been worked out in \cite[Sec.~3.2]{LeeuAH}. Again, for the sake of
completeness, we recall them here in a slightly different fashion.
Let us again consider the cell in
Figure~\ref{fig:1a}, filled by a non-negative integer $m$, and
labelled by the partitions $\rh$, $\mu$, $\nu$, where
$\rh\subseteq\mu$ and $\rh\subseteq\nu$, $\mu$ and $\rh$ differ by a
vertical strip, and $\nu$ and $\rh$ differ by a vertical strip.
Then, with the usual truth function $\chi(\mathcal A)=1$ if $\mathcal A$ is
true and $\chi(\mathcal A)=0$ otherwise, 
$\la$ is determined by the following algorithm:
\begin{enumerate}
\item[(F$^4$0)]Set $\text {CARRY}:=m$ and $i:=1$.
\item[(F$^4$1)]Set
$\la_i:=\max\{\mu_i,\nu_i\}+\min\{\chi(\rh_i=\mu_i=\nu_i),\text {CARRY}\}$.
\item[(F$^4$2)]If $\la_i=0$, then stop. The output of the algorithm
is $\la=(\la_1,\la_2,\dots,\la_{i-1})$. If not, then
set $\text {CARRY}:=\text {CARRY}-
\min\{\chi(\rh_i=\mu_i=\nu_i),\text {CARRY}\}+
\min\{\mu_i,\nu_i\}-\rh_i$ and $i:=i+1$.
Go to (F$^4$1).
\end{enumerate}

Conversely,
given $\mu,\nu,\la$, where $\mu\subseteq\la$ and $\nu\subseteq\la$,
where $\la$ and $\mu$ differ by a vertical strip, and 
where $\la$ and $\nu$ differ by a vertical strip, 
the backward algorithm works in the following way:
\begin{enumerate}
\item[(B$^4$0)]Set $i:=\max\{j:\la_j\text{ is positive}\}$ and $\text
{CARRY}:=0$.
\item[(B$^4$1)]Set $\rh_i:=\min\{\mu_i,\nu_i\}-
\min\{\chi(\mu_i=\nu_i=\la_i),\text {CARRY}\}$.
\item[(B$^4$2)]Set $\text {CARRY}:=\text {CARRY}-
\min\{\chi(\rh_i=\mu_i=\nu_i),\text {CARRY}\}+
\la_i-\max\{\mu_i,\nu_i\}$ and $i:=i-1$.
If $i=0$, then stop. The output of the algorithm is
$\rh=(\rh_1,\rh_2,\dots)$ and $m=\text {CARRY}$. If not, go to (B$^4$1).
\end{enumerate}

The corresponding extension of Theorem~\ref{thm:1} then reads as follows.

\begin{Theorem} \label{thm:1d}
Let $F$ be a Ferrers shape given by the $D$-$R$-sequence
$w=w_1w_2\dots w_k$.
Fillings of $F$ with non-negative integers 
are in bijection with sequences
$(\emptyset=\la^0,\la^1,\break\dots,\la^{k}=\emptyset)$, where
$\la^i/\la^{i-1}$ is a vertical strip if $w_i=R$, whereas
$\la^{i-1}/\la^i$ is a vertical strip if $w_i=D$.
\end{Theorem}

In the case that the Ferrers shape is a rectangle,
the bijection in Theorem~\ref{thm:1d} is equivalent to the
following insertion algorithm.
In order to describe this algorithm, 
we transform again our filling in
Figure~\ref{fig:9} into a two-rowed array. We do this in the
same way as in the insertion procedure modelling the third variant.
Namely, for an entry $m$ in the $i$-th row (from below) and $j$-th
column (from the left), we record $m$ pairs $\binom ji$. The pairs
are then ordered such that the entries in the top row are weakly increasing, 
in the bottom row, however, entries must be {\it decreasing}
below equal entries in the top row. Clearly, in the case of our
running example, this yields again the two-rowed array
\eqref{eq:2rc}.

Now we apply column insertion to the bottom entries to construct the
semistandard tableau $P^t$ and use the top entries to record
the insertions in the semistandard tableau $Q^t$. 
If we apply this procedure to our two-rowed array in \eqref{eq:2rc},
then we obtain the pair
\begin{equation} \label{eq:P,Q2}
(P^t,Q^t)=\left(\begin{matrix} 1\ 1\ 3\ 4\\2\ \hphantom{2\ 2\ 2}
\end{matrix}\ ,\
\begin{matrix} 1\ 1\ 1\ 2\\2\ \hphantom{2\ 2\ 2}
\end{matrix}\right),
\end{equation}
which is indeed in agreement with the increasing sequences
of partitions along the right side and the top side of the
rectangle in Figure~\ref{fig:9}, respectively, if one transposes both
arrays in \eqref{eq:P,Q2}.
(An equivalent insertion algorithm is described in
\cite[App.~A.4.1]{FultAC}, and is called ``Burge correspondence"
there since it appears as an aside in \cite{BurgAA}. 
The same terminology is used in \cite[Sec.~3.2]{LeeuAH}.)

\medskip
As before, the bijection in Theorem~\ref{thm:1d} has also a {\it global\/}
description, as a consequence of Greene's theorem.
In view of previous definitions, it should be clear what we mean
by {\it ne-chains} and {\it SE-chains}. We are then ready to
formulate the corresponding result.

\begin{Theorem} \label{thm:3d}
Given a diagram with empty partitions labelling all the
corners along the left side and the bottom side of the Ferrers shape, 
which has been completed according to dual RSK',
the partition $\la=(\la_1,\la_2,\dots,\la_\ell)$ labelling corner $c$ 
satisfies the following two properties:
\begin{enumerate}
\item [(G$^4$1)]Fix a positive integer $k$ and consider collections
of $k$ ne-chains with the property that no entry $e$ can be in more
than $e$ of these ne-chains.
Then the maximal cardinality of the multiset\/ {\em(!)} union of the ne-chains 
in such a collection of ne-chains
situated in the rectangular region to the left and
below of $c$ is equal to $\la_1+\la_2+\dots+\la_k$.
\item [(G$^4$2)]For any $k$, the maximal sum of all the entries
in a collection of $k$ pairwise disjoint
SE-chains situated in the rectangular region to the left and
below of $c$ is equal to $\la'_1+\la'_2+\dots+\la'_k$, where $\la'$
denotes the partition conjugate to $\la$.
\end{enumerate}
In particular, $\la_1$ is the length of the longest ne-chain
in the rectangular region to the left and below of $c$, and $\la'_1$
is the length of the longest SE-chain in the same
rectangular region.
\end{Theorem}

\subsection{Consequences on the enumeration of fillings with
restrictions on their increasing and decreasing chains}

We are now ready to state and prove the two extensions of
Theorem~\ref{thm:2} which result from Theorems~\ref{thm:3a},
\ref{thm:3b} and \ref{thm:3d}. In the statement,
in analogy to previous notation, 
we write $N^*(F;n;\NE=s,se=t)$ for the number
of fillings of the Ferrers shape $F$ with non-negative integers
with sum of entries equal to $n$
such that the longest
NE-chain has length $s$ and the longest se-chain, the
smallest rectangle containing the chain being contained in $F$, has
length $t$. The notation $N^*(F;n;ne=s,\SE=t)$ has the obvious
analogous meaning. 
Furthermore, we write $N^{01}(F;n;nE=s,\Se=t)$ for the number
of $0$-$1$-fillings of the Ferrers shape $F$ 
with sum of entries equal to $n$
such that the longest
nE-chain has length $s$ and the longest Se-chain, the
smallest rectangle containing the chain being contained in $F$, has
length $t$. The notation $N^{01}(F;n;\Ne=s,sE=t)$ has the obvious
analogous meaning.
We then have the following extensions of Theorem~\ref{thm:2}.

\begin{Theorem} \label{thm:2a}
For any Ferrers shape $F$ and positive integers $s$ and $t$, we have
\begin{equation} \label{eq:NES1}
N^*(F;n;\NE=s,se=t)=N^*(F;n;ne=t,\SE=s)
\end{equation}
and
\begin{equation} \label{eq:NES2}
N^{01}(F;n;nE=s,\Se=t)=N^{01}(F;n;\Ne=t,sE=s).
\end{equation}
\end{Theorem}

\begin{proof} 
We sketch the proof of \eqref{eq:NES1}.
We define a bijection between the fillings
counted by $N^*(F;n;\NE=s,se=t)$ and those counted by 
$N^*(F;n;ne=t,\SE=s)$.
Let the Ferrers shape $F$ be given by the $D$-$R$-sequence
$w=w_1w_2\dots w_k$.
Given a filling counted by $N^*(F;n;\NE=s,se=t)$ we apply the
mapping of the proof of Theorem~\ref{thm:1a}. Thus, we obtain a
sequence
$(\emptyset=\la^0,\la^1,\dots,\la^{k}=\emptyset)$, where
$\la^i/\la^{i-1}$ is a horizontal strip if $w_i=R$, whereas
$\la^{i-1}/\la^i$ is a horizontal strip if $w_i=D$. 
Since the sum of the entries in the filling was $n$, in the 
sequence the sum of the amounts of ``rises" $\la^{i-1}\subseteq\la^i$
(the amount is the number of squares by which $\la^i$ and $\la^{i-1}$
differ; equivalently, the sum of the amounts of ``falls" 
$\la^{i-1}\supseteq\la^i$ is exactly $n$).
Moreover, by Theorem~\ref{thm:3a}, we have $\la^i_1\le s$ and
$(\la^j)'_1\le t$ for all $i$ and $j$, with equality for at least one
$i$ and at least one $j$. Now we apply the inverse mapping from
Theorem~\ref{thm:1d} to
the sequence $(\emptyset=(\la^0)',(\la^1)',\dots,(\la^{k})'=\emptyset)$
of conjugate partitions. Thus, due to Theorem~\ref{thm:3d},
we obtain a filling counted
by $N^*(F;n;ne=t,\SE=s)$.

The proof of \eqref{eq:NES2} is completely analogous, but uses
the mapping from Theorem~\ref{thm:1b} to go from fillings counted by
$N^{01}(F;n;nE=s,\Se=t)$ to oscillating sequences of partitions,
and the inverse ``reflected" mapping from Section~\ref{sec:RSK'} 
to go from the sequence of
conjugate partitions to fillings counted by
$N^{01}(F;n;\Ne=s,sE=t)$. Theorem~\ref{thm:3b} is used to see how
the lengths of the nE-chains, Se-chains, Ne-chains, and sE-chains are
related to the partitions in the oscillating sequences of partitions.
\end{proof}

In the special case that $F$ is triangular,
it is straight-forward to use the bijections of this section
to formulate extensions of
Theorems~\ref{thm:4}--\ref{thm:6} to partitions of {\it multisets}.
We omit the details here for the sake of brevity. The reader will have no
difficulty to work them out.

\section{The big picture?}\label{sec:big}

In this paper we have been investigating which results on fillings,
where length restrictions are imposed on their chains,
can be obtained by Robinson--Schensted-like correspondences.
While I believe that the corresponding analysis here is (more or less) complete,
I believe at the same time that the obtained results constitute just a 
small section in a much larger field of phenomena that are concerned
with fillings which avoid certain patterns. The evidence that I can
put forward are recent results by 
Backelin, West and Xin \cite{BaWXAA},
Bousquet--M\'elou and Steingr\'\i msson \cite{BoStAA}
(generalizing previous results of Jaggard \cite{JaggAB}),
and Jonsson \cite{JonsAA}.
This, speculative, last section is devoted to a comparison of these
results with ours, and to posing some problems that suggest 
themselves in this context.

To begin with, in \cite{JonsAA} Jonsson considers $0$-$1$-fillings
of {\it stack polyominoes}. Here, a stack polyomino is the
concatenation of a Ferrers shape in French notation that has
been reflected in a vertical line with another (unreflected)
Ferrers shape in French notation. See Figure~\ref{fig:10}.a for an
example. Extending previously introduced notation,
we write $N^{01}(F;n;ne=s,se=t)$ for
the number of $0$-$1$-fillings of a stack polyomino 
$F$ with exactly $n$ $1$'s such that the longest
ne-chain, the smallest rectangle containing the chain being
contained in $F$, has length $s$ and the longest se-chain, the
smallest rectangle containing the chain being contained in $F$, has
length $t$.
Then Jonsson proves in \cite[Theorem~14]{JonsAA} (Theorem~13 of
\cite{JonsAA} contains even a refinement) that
\begin{equation} \label{eq:Jonsson1}
N^{01}(F;n;ne=s,se=*)=N^{01}(F';n;ne=s,se=*)
\end{equation}
if $n$ is maximal
so that $0$-$1$-configurations exist, the ne-chains of which are 
of length at most $s$, where $F'$ is the Ferrers shape which arises
by permuting the columns of $F$ so that they are ordered from the
longest column to the shortest. 
($se=*$ means that there is no restriction on
the length of the se-chains.)
Moreover, he reports that he and Welker
\cite{JoWeAA} extended this result to arbitrary $n$
(currently only for triangular shapes, but most likely an extension
to arbitrary stack polyominoes is possible without much difficulty).
Not only that, he also says that he expects these results
to remain true if $F$ is a moon polyomino, i.e., an arrangement of
cells such that along any row of cells and along any column of cells
there is no hole, and such that any two columns of $F$ have the
property that one column can be embedded in the other by applying a
horizontal shift. See Figure~\ref{fig:10}.b for an example.

\begin{figure}[h]
$$
\Einheit1cm
\Pfad(0,0),22122112121666161665555555\endPfad
\Pfad(0,1),1111111\endPfad
\Pfad(0,2),1111111\endPfad
\Pfad(1,3),11111\endPfad
\Pfad(1,4),1111\endPfad
\Pfad(3,5),11\endPfad
\Pfad(1,0),22\endPfad
\Pfad(2,0),2222\endPfad
\Pfad(3,0),2222\endPfad
\Pfad(4,0),22222\endPfad
\Pfad(5,0),222\endPfad
\Pfad(6,0),22\endPfad
\hbox{\hskip9cm}
\Pfad(0,1),221121221616616656555525\endPfad
\Pfad(1,1),1111\endPfad
\Pfad(0,2),111111\endPfad
\Pfad(0,3),111111\endPfad
\Pfad(2,4),111\endPfad
\Pfad(3,5),1\endPfad
\Pfad(1,1),22\endPfad
\Pfad(2,0),222\endPfad
\Pfad(3,0),2222\endPfad
\Pfad(4,0),22222\endPfad
\Pfad(5,1),222\endPfad
\hskip6cm
$$
\vskip10pt
\centerline{\small a. A stack polyomino
\hskip5cm
b. A moon polyomino}
\caption{}
\label{fig:10}
\end{figure}

Jonsson proves his result by an involved inductive argument which
does not shed any light why his result is true. Jonsson and Welker, on
the other hand, use machinery from commutative algebra to prove their
generalization. However, the most natural proof that one could think
of is, of course, a bijective one. This leads us to our first
problem.

\begin{Problem}
Find a bijective proof of \eqref{eq:Jonsson1}.
\end{Problem}

Where is the (at least, potential) connection to the material of our paper?
If we would apply Jonsson and Welker's results to 
a stack polyomino $F$ and its reflection in a vertical line, 
we would obtain that
\begin{equation} \label{eq:Jonsson2}
N^{01}(F;n;ne=s,se=*)=N^{01}(F;n;ne=*,se=s).
\end{equation}
Clearly, this is very close to the assertions
of Theorems~\ref{thm:2} and \ref{thm:2a}, and it
brings us to our next problem.

\begin{Problem} \label{prob:1}
Is it true that for any Ferrers shape {\em(}stack polyomino, moon
polyomino{\em)} $F$
and positive integers $s$ and $t$, we have
\begin{equation} \label{eq:Jonsson}
N^{01}(F;n;ne=s,se=t)=N^{01}(F;n;ne=t,se=s)\quad \text {\em?}
\end{equation}
\end{Problem}

The reader should observe that, although Theorems~\ref{thm:2} and
\ref{thm:2a} are very similar to \eqref{eq:Jonsson}, the bijections
used in their proofs cannot be used to prove \eqref{eq:Jonsson}.
Certainly, Theorem~\ref{thm:2} addresses $0$-$1$-fillings, but
$0$-$1$-fillings with the {\it additional condition} that in
each row and in each column there is at most one $1$. 
Also \eqref{eq:NES2} is about $0$-$1$-fillings, and indeed 
about $0$-$1$-fillings without
further restrictions. However, the
bijection does not keep track of lengths of ne-chains and se-chains,
but instead of lengths of nE-chains, Se-chains, Ne-chains and sE-chains.
On the other hand, \eqref{eq:NES1} does address ne-chains
and se-chains. Indeed, it would imply
\begin{equation} \label{eq:RSK3}
N^*(F;n;ne=s,se=*)=N^*(F;n;ne=*,se=s).
\end{equation}
As Jakob Jonsson (private communication) pointed out to me,
although this is not an assertion about $0$-$1$-fillings
but about {\it arbitrary} fillings, and although it does not imply
\eqref{eq:Jonsson}, it implies
at least \eqref{eq:Jonsson2}. Namely, consider the simplicial complex
$\De_{ne\le s}$ of $0$-$1$-fillings of $F$ with ne-chains being of
length at most $s$. (We refer the reader to \cite[Ch.~II]{StanAM} 
or \cite[Ch.~5]{BrHeAB}
for terminology and background on simplicial complexes and
Stanley--Reisner rings.)
Similarly, define the simplicial complex
$\De_{se\le s}$. There is an obvious bijection between monomials in
the Stanley--Reisner ring of $\De_{ne\le s}$ and {\it arbitrary}
fillings of $F$ with ne-chains being of length at most $s$, and 
there is an analogous, equally 
obvious bijection between monomials in
the Stanley--Reisner ring of $\De_{se\le s}$ and {\it arbitrary}
fillings of $F$ with se-chains being of length at most $s$.
Using this language, \eqref{eq:RSK3} says that there are as many
monomials of degree $n$ in the Stanley--Reisner ring of $\De_{ne\le
s}$ as there are monomials of degree $n$ in the Stanley--Reisner ring
of $\De_{se\le s}$. But this says that the Hilbert functions for
these two Stanley--Reisner rings are the same, which implies that the 
corresponding simplicial complexes have the same $h$-vector and
dimension, and 
hence the same $f$-vector. The latter means that,
for any positive integer $n$, there are as many
$0$-$1$-fillings with $n$ $1$'s and ne-chains being of length at most
$s$ as there are $0$-$1$-fillings with $n$ $1$'s and se-chains being
of length at most $s$, which is equivalent to \eqref{eq:Jonsson2}.

\begin{Problem}
Are there extensions of Theorems~{\em\ref{thm:2}} and {\em\ref{thm:2a}}
to stack polyominoes? To moon polyominoes? Are there analogues of
\eqref{eq:Jonsson1} for nE-chains, Se-chains, Ne-chains,
sE-chains, for arbitrary fillings?
\end{Problem}

The papers by 
Backelin, West and Xin \cite{BaWXAA}
and Bousquet--M\'elou and Steingr\'\i msson \cite{BoStAA}, on the other hand, consider
$0$-$1$-fillings of Ferrers shapes with the property that there is
{\it exactly} one $1$ in every row and in every column of the Ferrers shape.
As Mireille Bousquet--M\'elou pointed out to me, the bijection
behind Theorem~\ref{thm:2} provides in fact alternative proofs of the
results in these two papers. In the case of the main result in
\cite{BoStAA} this alternative proof is actually 
considerably simpler.

To be more specific, it is shown in \cite[Prop.~2.3]{BaWXAA} 
that the main
theorem in that paper, \cite[Theorem~2.1]{BaWXAA}, is implied
by the key assertion (see \cite[Prop.~2.2]{BaWXAA}) that 
there are as many such fillings of a Ferrers shape $F$ with longest
NE-chain of length $s$ as there are such fillings with longest SE-chain of
length $s$. This is indeed covered by Theorem~\ref{thm:2}, by
choosing $n$ to be equal to the maximum of width and height of
$F$. (Clearly, if width and height of a Ferrers shape $F$ are not
equal, then there are no $0$-$1$-fillings with exactly one $1$ in
every row and in every column of $F$.) 

Bousquet--M\'elou and Steingr\'\i msson \cite{BoStAA},
on the other hand, prove an
analogue of the main result of Backelin, West and Xin for {\it
involutions} (see \cite[Theorem~1]{BoStAA}). 
Their proof is based on a result due to Jaggard 
\cite[Theorem~3.1]{JaggAB} (see also \cite[Prop.~3]{BoStAA})
which shows that it is enough to prove that,
for any {\it symmetric} Ferrers shape $F$,
there are as many {\it symmetric} $0$-$1$-fillings 
of $F$ (as above, that is, with the property that there is exactly one $1$ in
every row and in every column of $F$) with longest
NE-chain of length $s$ as there are such symmetric 
$0$-$1$-fillings with longest SE-chain of
length $s$. (Here, ``symmetric" means that the filling, respectively
the Ferrers diagram, remains
invariant under a reflection with respect to the main diagonal.
In the case of the Ferrers diagram, a more
familiar term is ``self-conjugate.") 
It is now not too difficult to see that the growth diagram bijection
from Section~\ref{sec:growth} proving Theorems~\ref{thm:1}
and \ref{thm:2} provides a proof of the latter key assertion.
Namely, under the bijection in the proof of Theorem~\ref{thm:1}, {\it
symmetric} $0$-$1$-fillings correspond to {\it symmetric} sequences 
$(\emptyset=\la^0,\la^1,\dots,\la^{k}=\emptyset)$ as described in
the theorem, that is, to sequences with the property
$\la^i=\la^{k-i}$ for all $i$. Hence, the bijection in
Theorem~\ref{thm:2} restricts to a bijection between {\it symmetric}
$0$-$1$-fillings. 

The above arguments show that our bijection proves actually
a refined version of \cite[Prop.~2.2]{BaWXAA}, in which
one can keep track of the length of NE-chains and SE-chains
{\it at the same time}. More generally, these arguments prove
a version of Theorem~\ref{thm:2} for {\it symmetric}
$0$-$1$-fillings. Namely, for a symmetric (= self-conjugate) 
Ferrers diagram $F$
let us write $N_{\text{sym}}(F;n;\NE=s,\SE=t)$ for the number
of {\it symmetric} 
$0$-$1$-fillings of the Ferrers shape $F$ with exactly $n$ $1$'s,
such that there is at most one $1$ in each column and in each row,
and such that the longest
NE-chain has length $s$ and the longest SE-chain, the
smallest rectangle containing the chain being contained in $F$, has
length $t$. Then we have the following theorem.

\begin{Theorem} \label{thm:2sym}
For any symmetric
Ferrers shape $F$ and positive integers $s$ and $t$, we have
$$N_{\text{sym}}(F;n;\NE=s,\SE=t)=N_{\text{sym}}(F;n;\NE=t,\SE=s).$$
\end{Theorem}

Although the above arguments treat the results from \cite{BaWXAA} and
\cite{BoStAA} in a uniform manner, and although it seems that
they elucidate them to a consideryble extent, 
there remain still some mysteries
that seem worth to be investigated further. These mysteries concern
the transformation due to Backelin, West and Xin on which the proofs
in \cite{BaWXAA} and \cite{BoStAA} are based. In order to explain
this, let us recall that neither in
\cite{BaWXAA} nor in \cite{BoStAA} are the above enumerative key
assertions proved by a direct bijection, but in a round-about way
which, in the end, is again based on a bijection,
but a bijection between different sets of objects.
This bijection is described in Section~3 of \cite{BaWXAA} and is
called ``transformation $\phi$" in \cite{BoStAA}. Let us call it the
``BWX transformation." As is explained in \cite[Prop.~2.4]{BaWXAA},
the BWX transformation induces a bijection for the first of the two
above key assertions, and vice versa, any bijection for the key
assertion induces a bijection that could replace the BWX
transformation. Thus, the question arises whether the BWX
transformation and our growth diagram bijection have something to do
which each other, via the link in \cite[Prop.~2.4]{BaWXAA}.
This question is of particular interest in connection with
the article \cite{BoStAA}: the key result in the latter paper is the
proof that the BWX transformation commutes with reflection.
This turns out to be extremely hard. On the other hand, that
the growth diagram bijection in Theorem~\ref{thm:2} commutes
with reflection is a triviality.

\begin{Problem}
What is the relation between the BWX transformation and
the growth diagram bijections in Section~\ref{sec:growth}?
\end{Problem}

Interestingly, in a not at all equally precise, but somewhat
speculative form, 
this question was already posed at the end of the Introduction
of \cite{BoStAA}. Now, in the light of the above remarks
and observations,
this speculation seems to get firm ground.

For the convenience of the reader, we also state the ``symmetric"
version of Theorem~\ref{thm:2a} explicitly, which follows in a
manner completely analogous to the way 
Theorem~\ref{thm:2sym} followed from the bijection in
Theorem~\ref{thm:2}.
In the statement, $N_{\text{sym}}^*(F;n;\NE=s,se=t)$ denotes 
the number
of {\it symmetric} fillings of the symmetric Ferrers shape 
$F$ with non-negative integers
with sum of entries equal to $n$
such that the longest
NE-chain has length $s$ and the longest se-chain, the
smallest rectangle containing the chain being contained in $F$, has
length $t$. The notation 
$N_{\text{sym}}^*(F;n;ne=s,\SE=t)$ has the obvious
analogous meaning. 
Furthermore, $N_{\text{sym}}^{01}(F;n;nE=s,\Se=t)$ 
denotes the number
of {\it symmetric} $0$-$1$-fillings of the symmetric Ferrers shape $F$ 
with sum of entries equal to $n$
such that the longest
nE-chain has length $s$ and the longest Se-chain, the
smallest rectangle containing the chain being contained in $F$, has
length $t$. The notation $N_{\text{sym}}^{01}(F;n;\Ne=s,sE=t)$ has the obvious
analogous meaning.

\begin{Theorem} \label{thm:2asym}
For any symmetric Ferrers shape $F$ and positive integers $s$ and $t$, we have
\begin{equation} \label{eq:NES1sym}
N_{\text{sym}}^*(F;n;\NE=s,se=t)=N_{\text{sym}}^*(F;n;ne=t,\SE=s)
\end{equation}
and
\begin{equation} \label{eq:NES2sym}
N_{\text{sym}}^{01}(F;n;nE=s,\Se=t)=N_{\text{sym}}^{01}(F;n;\Ne=t,sE=s).
\end{equation}
\end{Theorem}

If one goes back to the papers by
Backelin, West and Xin, and by Bousquet--M\'elou and
Steingr\'\i msson, which are about permutations with forbidden
patterns, then one may wonder if there are also
theorems in the case of arbitrary fillings for more
general patterns than just increasing or decreasing
patterns, or in the case
of $0$-$1$-fillings
where one relaxes the condition that every row and column must
contain exactly one $1$. 
On the other hand,
if one recalls the results by Jonsson and Welker, then one is
tempted to ask whether the above theorems may be extended,
in one form or another, to more general shapes.
Thus, we are led to the following problem.

\begin{Problem}
Are there extensions of Theorems~{\em\ref{thm:2}}, 
{\em\ref{thm:2a}--\ref{thm:2asym}},
or of {\em\eqref{eq:Jonsson}} to more general patterns? To more
general shapes? 
\end{Problem}

Finally one may also ask the question whether there is anything
special with $0$-$1$-fillings as opposed to arbitrary fillings, or
whether one can also obtain results in the spirit of this paper for
fillings where the size of the entries is at most $m$, for an {\it
arbitrary} fixed $m$.

\begin{Problem}
Can one extend the results for $0$-$1$-fillings to
fillings with entries from $\{0,1,2,\dots,m\}$?
\end{Problem}

Further papers which should be considered in the present context are 
\cite{ChMYAA,KlazAD,KlazAE}.

\section*{Acknowledgments}
The idea for this article arose in discussions with Donald
Knuth at the Institut Mittag--Leffler during the
``Algebraic Combinatorics" programme in Spring 2005,
which are herewith gratefully acknowledged.
I would like to thank Anders Bj\"orner and Richard
Stanley, and the Institut Mittag--Leffler, 
for giving me the opportunity
to work in a relaxed and inspiring atmosphere
in the framework of this programme.
Furthermore, I am indebted to Mireille Bousquet--M\'elou, Jakob Jonsson and
Marc van Leeuwen for many useful suggestions for improvements, 
corrections, and for
pointing out short-comings in a first version, which helped to
improve the contents of this paper considerably.

\end{document}

